# La Cosmología y los matemáticos

por

**José M. M. Senovilla**

Se presentan algunos de los hitos históricamente relevantes y llevados a cabo o instigados, de manera esencial, por matemáticos en la creación, avance y desarrollo de la cosmología como disciplina científica. Asimismo, se detalla la estrecha relación entre las matemáticas y la cosmología, a través de la geometrización de ésta llevada a cabo por Einstein con su teoría de la relatividad general y la colaboración posterior de ilustres matemáticos del siglo XX.

## 1 Introducción

Como es notorio, durante 2005 celebramos el Año Mundial de la Física con motivo del centenario del *annus mirabilis* del grandísimo Albert Einstein (Ulm, Alemania, 1879 – Princeton, EE UU, 1955). En 1905 Einstein, a la sazón empleado de una oficina de patentes en Berna, publicó cuatro artículos revolucionarios, visionarios, bellos, de consecuencias devastadoras para la ciencia –¡y las matemáticas!– moderna.

La lista de revoluciones científicas provocadas, comenzadas o catalizadas por Einstein es tan extensa que seguramente resulta inverosimil. Probablemente por ello es el más famoso de los poquísimos científicos conocidos por los legos. La mayoría de dichas revoluciones se fraguaron antes de o durante el celebrado 1905: los fundamentos de las teorías cinética y estadística; la explicación del movimiento browniano de partículas en suspensión de un líquido; el efecto fotoeléctrico y la introducción de los *quanta* de energía, lo que llevó a la teoría *cuántica* de la radiación y posteriormente de la materia; la teoría de la relatividad especial o restringida. Hay no obstante otras revoluciones einsteinianas que tuvieron lugar bastante después, las más conocidas de éstas son, naturalmente, la teoría de la relatividad general y las bases de la estadística cuántica.

Por si el lector no está abrumado aún, permítanme la originalidad de hablar aquí, empero, de *otra* revolución fundamental iniciada por Einstein, con consecuencias del mismo calado y probablemente mayor alcance, que trastornó completamente nuestra manera de conocer y analizar el Cosmos: el nacimiento de la Cosmología como una genuina disciplina científica.

Teniendo en cuenta que este año se celebra también el cincuentenario de la desaparición tanto de Einstein como del influyente matemático, y colega suyo en Princeton, Hermann Weyl (Elsmhorn, Alemania, 1885 – Zürich, Suiza, 1955), quien contribuyó también al desarrollo de la relatividad general y la geometría del espacio-tiempo, este artículo se adentrará en la relación entre la cosmología y las matemáticas, y en las contribuciones, absolutamente imprescindibles, de grandes matemáticos a aquella. Sirva como homenaje a Einstein, Weyl, Friedman (San Petersburgo 1888 – Leningrado 1925) de quién ahora se cumplen 80 años de su muerte, y a los demás matemáticos que irán apareciendo en el texto de este artículo: Pitágoras, Ptolomeo, Newton, Bianchi, Riemann, Minkowski, Grossmann, Christoffel, Levi-Civita, Ricci-Curbastro, Jacobi, Lanczos, Robertson, Walker, Gödel, Taub...

Lo que sigue contiene exclusivamente hechos conocidos, pero teñidos de un punto de vista personal que me permitirá resaltar los pasos fundamentales y los acontecimientos extraordinarios que rodearon esta apasionante historia.

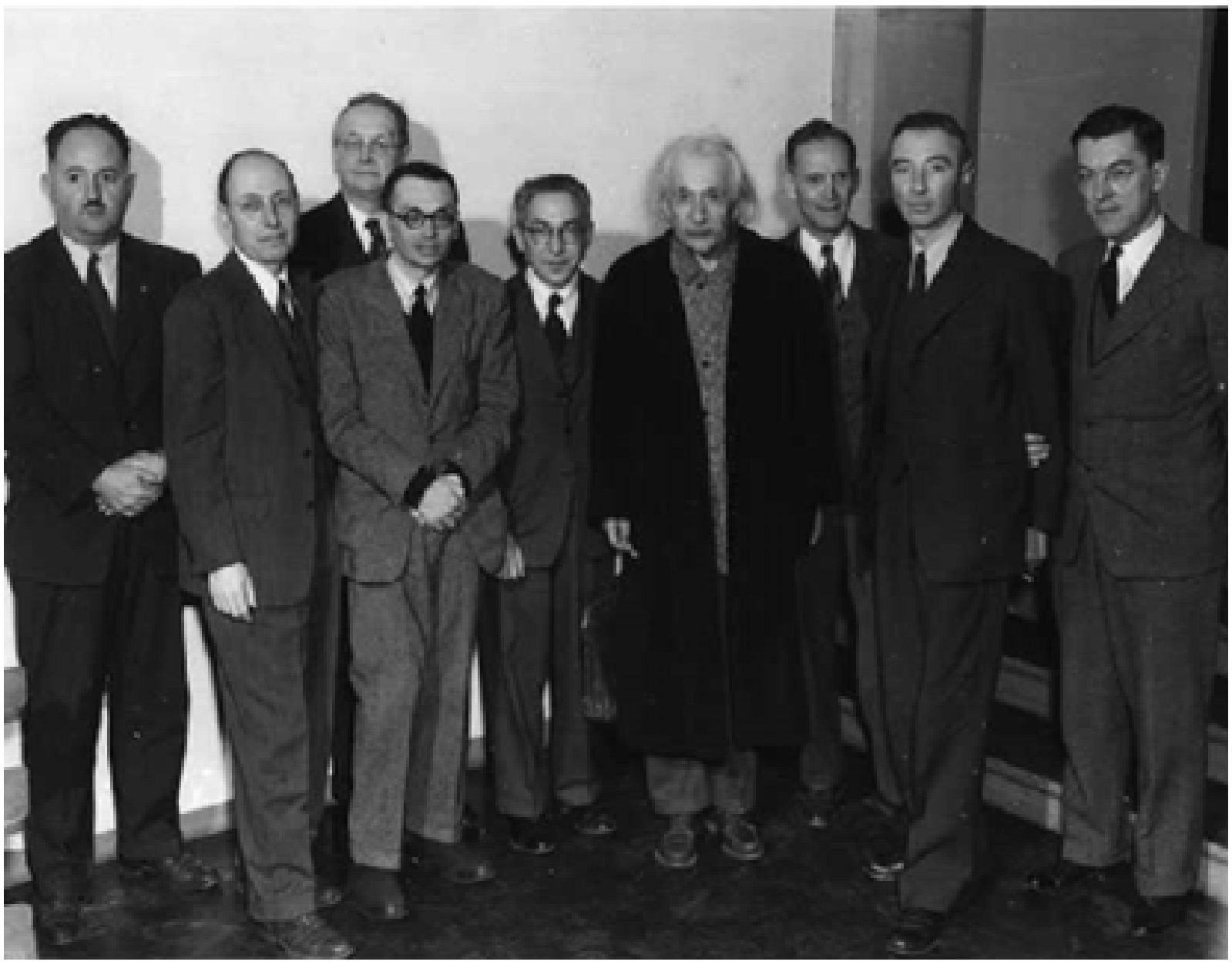

Celebración del septagésimo cumpleaños de Einstein en Pasadena
(de izq. a dch.) H.P. Robertson, E. Wigner, H. Weyl, K. Goedel, I.I. Rabi, A. Einstein, R. Ladenburg, J.R. Oppenheimer y G.M. Clemence.

## 2 Antes de 1905

La historia de la ciencia (y de la filosofía) está salpicada de visiones del cosmos, intuiciones sobre su forma, sus propiedades y composición. Hay por ejemplo innumerables *cosmogonías* antiguas, mayormente basadas en creencias o razones de fe, cuyo principal objetivo era servir de ambiente para los fenómenos, identificar el Universo como marco en el que ocurren las cosas. Se pueden citar, como casos destacados

- las cosmogonías de la escuela naturalista de Mileto entre 560 y 525 a.c. –como las de Tales, Anaximandro o Anaxímenes–, en las que la Tierra tenía forma de disco plano o de rodaja de sandía, y los cuatro elementos se distribuían debajo de una "esfera celeste";

- la cosmogonía de la escuela de Elea, debida a Parménides *circa* 530 a.c., en la que la Tierra era redonda preservando las demás propiedades de las cosmogonías previas;

- la del pitagorismo, escuela de resonancia matemática, atribuída a Filolao, en la que los números tuvieron una importancia fundamental ya que para mantener la perfección del número 10 se creó un fantasmagórico planeta llamado *Antichton* o anti-Tierra que resguardaba la Tierra del fuego central. Era éste, por tanto, un modelo no geocéntrico, aunque tampoco heliocéntrico, ya que el Sol giraba también en torno al mencionado fuego;

- el primer sistema heliocéntrico se debe al gran Aristarco de Samos, quien *circa* 260 a.c. propuso que la Tierra era esférica, giraba en torno al Sol y también sobre sí misma, e incluso puso el eje de rotación de la Tierra inclinado respecto del plano de rotación en torno al Sol (hoy eclíptica).

Como vemos, las ideas de Aristarco de Samos se adelantaron muchos siglos a las de Copérnico, aunque se ha de señalar que tuvieron un impacto nimio debido a su "evidente sinrazón".

Aparte de lo anterior, se han de mencionar también los pocos marcos teóricos surgidos antes de los siglos XIX-XX para explicar el Cosmos. Seguramente el primero de relevancia fue el sistema ideado por Claudio Ptolomeo (Hermiou, Egipto, 85—Alejandría 165) en el famoso *Almagesto*, donde se da un tratado detalladísimo y completo de la "teoría de los cielos", como él mismo la llamó, fundando la astronomía, que permitió predecir el movimiento de los cuerpos celestes conocidos y llegar a una mejor comprensión del mundo. Como se sabe, todo ello basado en un modelo geocéntrico, de éxito casi sin par (probablemente sólo los *Elementos* de Euclides hayan estado en uso por un periodo mayor de tiempo), y que llegó prácticamente intacto al siglo XVI.

El siguiente marco teórico, y probablemente el de mayor relevancia histórica, surge después de la revolución heliocéntrica llevada a cabo por Copérnico (Torun, Polonia 1473 – Frombork 1543), continuada, defendida e impulsada por el nacimiento del método científico y los telescopios con Galileo (Pisa, 1564 – Arcetri 1642), ayudada por las precisas observaciones de Tycho Brahe (Knudstrup, Dinamarca, 1546 – Praga, Bohemia, 1601) y finalmente asentada con las magníficas leyes de Kepler (Weil der Stadt, Württemberg, 1571 – Regensburg, Alemania, 1630) del movimiento planetario. Nos referimos, claro está, al marco teórico creado en la genial síntesis desarrollada por Isaac Newton (Woolsthorpe, Inglaterra, 1643 – Londres 1727) en sus *Principia*, usando sus técnicas matemáticas del cálculo de fluxiones (hoy cálculo diferencial). Si el cálculo diferencial fue un paso gigantesco en el desarrollo de las matemáticas, no hay quizás adjetivo suficiente para resaltar lo que significaron los *Principia* en la historia de la física. Resaltemos aquí que las leyes de Kepler se *deducen* dentro del marco newtoniano basándose simplemente en el cálculo diferencial y en la famosa ley newtoniana de la gravitación universal,

$$\vec{F} = \frac{GMm}{r^2} \frac{\vec{r}}{r}, \tag{1}$$

que cuantifica la fuerza de atracción entre dos cuerpos cualesquiera con masas $M$ y $m$, siendo $G$ la constante de la gravitación universal y $\vec{r}$ el radio vector que une (los centros de masas de) los dos cuerpos. Los movimientos astronómicos quedaban así "explicados" dentro de la física clásica de Newton.

Dicho lo anterior, hay que resaltar empero que el marco teórico newtoniano presuponía, como hipótesis esencial (basada supuestamente en la experiencia) la existencia de un espacio y un tiempo absolutos, inmutables e inalterables, que no estaban influidos por ningún agente material, y eran infinitos. Este "aburrido" espacio-tiempo newtoniano estaba dado *a priori*, y servía de campo de referencia para el movimiento y la existencia de los objetos reales. Por ello mismo, el cosmos como un todo no se estudiaba, o sea, el Universo como tal no era *objeto* de estudio científico.

Es necesario reseñar ahora dos pasos menos conocidos pero de gran importancia para la historia que nos ocupa. El primero se debió a Richard Bentley (Oulton, Inglaterra 1662 – Cambridge 1742), contemporáneo de Newton y un hombre por lo demás disparatado, cuyo afán era probar la existencia de Dios usando la física newtoniana, quien puso de manifiesto que un Universo finito es inestable según la ley de Newton (1). Más aún, la misma conclusión se sigue si, a pesar de que el espacio fuera infinito, toda su materia estuviera concentrada sólo en una región finita, ya que toda la materia colapsaría por efecto de la atracción gravitatoria. Por esta razón el propio Newton sugirió que la distribución de la materia cósmica tenía que estar distribuida por todas partes, ser infinita y, probablemente entonces, estable de acuerdo con su ley (1), de manera que la fuerza neta resultante de todo el Universo sobre la Tierra pudiera considerarse como únicamente la del Sol. Como veremos después, esta conjetura resultó ser completamente insostenible.

El segundo paso que es necesario recordar, una de las ideas más lúcidas y visionarias en lo que concierne a la organización del cosmos, se debió al grandísimo filósofo Immanuel Kant (Königsberg, Alemania, 1724 – 1804), que en 1755 publicó en su tratado *Historia natural universal y teoría de los cielos* la idea de que, igual que el Sol formaba parte de un sistema aislado de estrellas con forma de disco, era probable que las nebulosas elípticas que podían observarse en el firmamento fuesen agrupaciones estelares del mismo tipo. Kant llamó a estas nebulosas "universos isla" haciendo una sorprendente predicción que se confirmaría ¡170 años más tarde!

A caballo entre los siglos XIX y XX aparecieron las primeras voces con críticas científicas de calado que abrirían la puerta al posterior desarrollo de la Cosmología. Una crítica de particular relevancia fue la realizada por Ernst Mach (Chirlitz-Turas, Austria 1838 – Vaterstetten, Alemania 1916), quien con demoledor estilo desmenuzó el problema que representaba el concepto de la inercia, definida *respecto del espacio absoluto*, en la teoria newtoniana. Mach discutió brillante y razonadamente que la inercia de un cuerpo no debía ser referida al espacio absoluto sino que, de hecho, ha de ser una propiedad relativa con *respecto de la distribución total de materia en el resto del Universo*. Esta idea influyó decisivamente en Einstein y sus primeras concepciones acerca del cosmos.

Otra crítica de similar importancia fue la denominada *paradoja de Seeliger*, que no es otra cosa que la concreción, con argumentos cuantitativos decisivos, de la inestabilidad del cosmos newtoniano señalada por Bentley y mencionada más arriba. Si hoy en día es evidente que la gravedad es la fuerza más relevante para el estudio del Universo, lo era más aún en el siglo XIX y principios del XX cuando sólo se conocían dos interacciones en la naturaleza, la electromagnética y la gravitatoria. Era por ello imprescindible usar una teoría de la gravitación en cualquier estudio del cosmos como un todo. La teoría vigente era la de la gravitación universal de Newton, pero conducía a conclusiones insostenibles como las mencionadas. La conjetura de Newton para refutar los argumentos de Bentley resultó ser totalmente falsa, como demostró rigurosamente Hugo von Seeliger (Biala bei Bielitz, Austria 1849 – Munich, Alemania 1924) en 1894: si la ley (1) valiera para el cosmos, el campo gravitatorio sería infinito por doquier. Curiosamente, Seeliger trató de modificar la ley (1) de la manera siguiente

$$\vec{F} = \frac{GMm}{r^2} e^{-kr} \frac{\vec{r}}{r} \tag{2}$$

para una constante positiva $k$. Esta nueva ley todavía permite explicar el comportamiento de los planetas en el sistema solar, pero al mismo tiempo la gravedad deviene muy débil a grandes distancias y es factible que distribuciones materiales a gran escala sean estables, incluso en el caso de distribuciones no homogéneas. Desafortunadamente, la modificación (2) de la ley de la gravitación no estaba exenta de problemas, como la inmensa mayoría de las modificaciones *ad hoc* con parámetros ajustables de cualquier teoría física. Hay que enfatizar que la paradoja de Seeliger y su propuesta solución no conllevaba una búsqueda de la forma o de las propiedades del Universo, no era un intento de construir una teoría cosmológica; simplemente sirvió como crítica a la visión newtoniana en la que se daba por supuesto que el espacio era infinito.

Llegados a este punto, hacía falta una nueva teoría de la gravitación para pavimentar el camino hacia una nueva Cosmología con estatus de ciencia. Esto es precisamente lo que creó Einstein con su relatividad general, pero no con la intención de resolver los problemas mencionados, sino por razones de otra índole, de carácter primordialmente estético, de simetría y completitud. Claro está que si, de paso, se podían resolver la paradoja de Seeliger y la problemática de la inercia de Mach, miel sobre hojuelas.

## 3 1905-1916: relatividad general y geometría lorentziana subyacente

Especialmente este año sabemos que, en 1905, Einstein llevó a cabo una de sus más famosas revoluciones por medio de la teoría de la relatividad especial, en la que se mezclaron el espacio y el tiempo y se tendió un puente que enlazaba la electrodinámica con la mecánica, dando carácter absoluto al valor de la velocidad de la luz en el vacío y elevando el principio de relatividad entre sistemas de referencia inerciales a la categoría de axioma indiscutible. Las consecuencias de este avance fueron, son y serán trascendentales para toda la física...

¡Y también para las matemáticas! Para empezar a saborear las decisivas e ingentes contribuciones de los matemáticos al desarrollo de la relatividad y, *a fortiori*, de la cosmología, nos encontramos con la unión en una única entidad –el espacio-tiempo– que realizó Hermann Minkowski (Alexotas [Kaunas], Rusia 1864 – Göttingen, Alemania 1909) en 1907. Este continuo espacio-tiempo seguía teniendo un carácter absoluto, e infinito, pero el espacio y el tiempo por separado dejaban de tener ese sentido y quedaban indisolublemente unidos para siempre (o eso creemos ahora). En sus propias palabras, de su más famosa cita:

> *A partir de ahora, el espacio por sí solo, y el tiempo por sí mismo, se han desvanecido en meras sombras y solamente una especie de mezcla de los dos existe por derecho propio.*

Como se sabe, Minkowski propuso entender la relatividad especial y la electrodinámica desarrollada por Lorentz, Poincaré y Einstein utilizando un continuo espacio-temporal. Usando coordenadas cartesianas $\{x, y, z\}$ para describir el espacio tridimensional y $t$ para denotar el tiempo en un sistema de referencia inercial[1], Minkowski propuso considerar un *elemento de línea* dado

[1] En física se puede definir teóricamente un sistema de referencia inercial como aquel en el que los cuerpos abandonados a sí mismos, sin agentes externos obrando sobre ellos, están en

por

$$ds^2 = -c^2dt^2 + dx^2 + dy^2 + dz^2\,, \tag{3}$$

siendo $c$ una constante que tiene el valor de la velocidad de la luz en un tal sistema de referencia y en el vacío. Como se sabe, este valor de $c$ es una constante universal según la teoría relativista. El elemento de línea anterior mide el *intervalo* infinitesimal entre cualesquiera dos sucesos (dos puntos) con coordenadas $(t, x, y, z)$ y $(t + dt, x + dx, y + dy, z + dz)$. Como se ve, el intervalo "al cuadrado" $ds^2$ no es necesariamente positivo, debido al signo menos delante del primer término del miembro derecho de (3), y puede ser negativo o nulo. El primer caso corresponde a puntos separados *temporalmente* y mide el tiempo transcurrido (por ejemplo, entre Ud. mismo al leer *esto* y *esto otro*); el segundo caso corresponde a puntos separados de manera *luminosa*, o sea por un rayo de luz, de manera que ningún otro agente físico puede conectarlos (la velocidad de la luz es la máxima posible); finalmente, el caso tradicional donde el intervalo $ds^2$ es positivo corresponde a una separación *espacial* y mide la distancia tradicional entre dos puntos en el sistema de referencia dado.

Una manera más matemática de decir lo mismo es que el marco de la relatividad especial es una *variedad pseudo-riemanniana de tipo lorentziano*, $(\mathbb{R}^4, g)$ cuyo campo tensorial métrico $g$ adopta la forma $g = \mathrm{diag}(-1, 1, 1, 1)$ en la cobase natural $\{cdt, dx, dy, dz\}$. O mejor aún, el campo tensorial métrico es

$$g = -c^2dt \otimes dt + dx \otimes dx + dy \otimes dy + dz \otimes dz\,. \tag{4}$$

Así pues, el tiempo alcanza la misma categoría que el resto de coordenadas espaciales, pero con un *status* especial, el que le da el signo menos en la signatura de la métrica $g$, el signo menos en la forma diagonal de $g$. Obviamente, con una tal métrica, que no es definida positiva, existen vectores tangentes con módulo negativo, nulo o positivo, y *a posteriori* curvas diferenciables con vectores tangentes que pueden ser de estos tres tipos, llamadas temporales, luminosas y espaciales, respectivamente.

### 3.1 La idea clave: el principio de equivalencia

Todo lo anterior permitía simplificar la matemática subyacente a la relatividad y a su primogenitora, la electrodinámica, de manera que la mecánica y el electromagnetismo tenían un nuevo y más potente lenguaje natural. Pero, ¿qué pasaba con la Gravitación? ¿Y dónde quedaban los sistemas de referencia no inerciales[2]? ¿Como "ciudadanos" de segunda categoría? La respuesta reposo o en movimiento rectilíneo y uniforme; la realización material, o incluso la existencia real, de tales sistemas ideales sólo se puede conseguir o asegurar de forma aproximada.

[2]Un sistema no inercial es cualquiera que esté acelerado con respecto a uno inercial, por ejemplo el de un tiovivo, o de un cohete en su lanzamiento, o el de un tren frenando, arrancando o dando una curva.

completa a estas cuestiones básicas sólo alcanzaron la condición de verdadera teoría físico-matemática a finales de 1915, gracias, en particular, a la ayuda que Einstein recibió de su amigo y compañero de estudios, el matemático Marcel Grossmann (Budapest, Hungría 1878 – Zurich, Suiza 1936). De esto hablaremos brevemente en lo que sigue, pero antes de eso hay que resaltar que la intuición fundamental, el resultado esencial sobre el que posteriormente se construiría todo el edificio, el concepto básico detrás de la teoría de la Relatividad General, ésta la tuvo Einstein muy temprano: en 1907. Se llama el *Principio de Equivalencia*, y es una genial pero sencillísima idea que *interpreta* –o sea, traduce, explica el sentido de– la igualdad de las masas gravitatoria e inercial[3] (un hecho que Newton ya conocía), dotando de paso a los sistemas no inerciales de la misma respetabilidad que cualquier otro sistema, e incluyendo asimismo la gravitación dentro de una teoría relativista. En palabras del propio Einstein: "glücklichste Gedanke meines Lebens" ("la idea más lúcida de toda mi vida" –lo cual, tratándose de Albert Einstein, ¡es mucho!).

El principio de equivalencia es un principio físico que se puede expresar como:

> *cualquier sistema de referencia no inercial es* localmente *equivalente a un sistema inercial con un cierto campo de gravitación.*

Viceversa,

> *todo campo gravitatorio puede hacerse desaparecer* localmente *escogiendo un sistema de referencia adecuado*

–por ejemplo, poniéndose en caída libre en la Tierra–. Este magnífico principio, comprobado experimentalmente a día de hoy hasta la saciedad, se basa en la propiedad fundamental de los campos de gravedad, ya comprobada por Galileo, de que todos los cuerpos libres adquieren la misma aceleración independientemente de su masa y constitución (en el vacío). Una propiedad que es obviamente cierta también en los sistemas acelerados. Es evidente que el principio de equivalencia estaba implícito en la mencionada igualdad de las masas gravitatoria e inercial, porque de la segunda ley de Newton y la fórmula (1) se deduce:

$$m\vec{a} = \frac{GMm}{r^2}\,\frac{\vec{r}}{r} \qquad \Longrightarrow \qquad \vec{a} = \frac{GM}{r^2}\,\frac{\vec{r}}{r}\,.$$

Las masas inercial y gravitatoria se han cancelado a ambos lados de la ecuación debido a su igualdad, de donde finalmente se sigue que el campo vectorial de aceleración $\vec{a}$ es independiente de la masa $m$ del cuerpo sometido a la fuerza de la gravedad creado por la masa $M$. Siendo esto conocido durante muchísimos años (desde mediados del siglo XVII), y a pesar de que grandes científicos avisaron de la relevancia y de la necesidad de explicar este hecho[4], durante todo ese tiempo nadie supo extraer la información ni las consecuencias que de él se derivaban hasta la llegada de Einstein y su enunciado del principio de equivalencia. Seguramente esta es la marca de los mayores genios frente a los, "simplemente", grandísimos científicos.

[3] La masa gravitatoria de un cuerpo es la que aparece en la fórmula (1) y mide la intensidad con que un cuerpo reacciona a la atracción de la gravedad. La masa inercial mide la inercia de un cuerpo, su resistencia a cambiar de movimiento, y es la que aparece en la segunda ley de Newton: $\vec{F} = m\vec{a}$, donde $\vec{a}$ es la aceleración del cuerpo respecto del espacio absoluto.

### 3.2 La geometrización

El principio de equivalencia lleva casi inexorablemente a la *geometrización* de la gravedad: dado que todos los cuerpos se mueven igual independientemente de su masa, podemos imaginar que de hecho siguen trayectorias determinadas por "la forma" del espacio, que simplemente van por los caminos menos esforzados de un espacio que no es plano. Esto es análogo a lo que sucede, por poner un ejemplo, cuando se quiere ir óptimamente de un punto a otro de la superficie terrestre: para ello hay que seguir las geodésicas, que como se sabe son trozos de circunferencias de radio máximo. Hay que observar que esta geometrización no es posible para el electromagnetismo, porque el movimiento de las partículas cargadas sí que depende de su relación carga/masa. Diferentes cargas siguen por tanto diferentes caminos.

Veamos ahora cómo puede usarse el principio de equivalencia para deducir, en primera aproximación, el marco matemático de la nueva teoría de la relatividad general. Supongamos que queremos describir el sistema de referencia de un tiovivo que se mueve con velocidad angular constante $\omega$. Si suponemos que los caballitos giran con respecto a un sistema de referencia inercial tradicional cuyas coordenadas cartesianas son $\{t, x, y, z\}$, basta para ello con realizar el siguiente cambio elemental de coordenadas

$$T = t, \qquad \rho = \sqrt{x^2 + y^2}, \qquad \phi = \arctan \frac{y}{x} + \omega t, \qquad Z = z,$$

válido en un dominio adecuado. Entonces, el campo métrico (4) se transforma en

$$\begin{aligned} g = -(c^2 - \omega^2\rho^2) dT \otimes dT - \omega\rho^2(dT \otimes d\phi + d\phi \otimes dT) + \\ + d\rho \otimes d\rho + \rho^2 d\phi \otimes d\phi + dZ \otimes dZ, \end{aligned} \tag{5}$$

[4]Por ejemplo, el gran Hertz, en su obra *Über der Constitution der Materie* (1884), escribió: "ciertamente tenemos frente a nosotros dos propiedades totalmente fundamentales de la materia que deben considerarse como completamente independientes entre sí, pero en nuestra experiencia, y sólo experimentalmente, se nos aparecen como iguales. Esta correspondencia debe significar mucho más que ser un simple milagro ( . . . ) Debemos percatarnos claramente de que la correspondencia entre masa e inercia tiene que tener una explicación más profunda, y no puede despacharse como algo de poca importancia...". A esto se llama profunda intuición y buen olfato; a lo de Einstein, genialidad.

que mantiene signatura lorentziana correcta en $c > \omega\rho$. Como se ve, el campo tensorial $g$ deja de ser diagonal, sus componentes ya no son constantes y, de mayor interés físico, el dominio de validez del sistema coordenado pasa a ser finito: de hecho, la región $\omega\rho \to c$ es un *horizonte*, un concepto de fundamental importancia en relatividad. Usando el principio de equivalencia, se sigue por lo tanto que este tipo de propiedades de $g$ aparecerán en las regiones con un campo gravitatorio. En consecuencia, se postula que un campo gravitatorio vendrá determinado por un cierto campo tensorial $g$ que es arbitrario excepto por dos requisitos fundamentales: que sea no degenerado, $\det g \neq 0$, por lo que $g$ no tiene radical lo que asegura la cuadridimensionalidad efectiva y real del espacio-tiempo; y que $g$ tenga signatura lorentziana –o sea, hay un tiempo y tres dimensiones espaciales–.

Claro está, para tratar de describir la física en un sistema no inercial tal como el anterior, o en un campo gravitatorio real, hace falta saber cómo se transforman las leyes de la física. Además, hay que incorporar en el marco matemático la idea de localidad, o en otras palabras, el hecho de que en un campo gravitatorio real se puedan anular sus efectos localmente. Para lo primero, que se ha venido en llamar "covariancia general", implicando que las leyes de la física se han de saber escribir en cualquier base (y sistema de coordenadas) hacía falta usar el cálculo tensorial, o diferencial absoluto, desarrollado por Christoffel, Ricci-Curbastro y Levi-Civita. Para lo segundo, era necesaria la teoría general de los espacios curvos de Riemann. En ambos aspectos la colaboración de Grossmann con Einstein fue vital, razón por la cual la comunidad relativista le sigue rindiendo homenaje actualmente mediante un congreso trianual que lleva su nombre.

Antes de seguir adelante, hay que hacer una pequeña pero decisiva puntualización acerca del uso de la palabra *local* en los enunciados precedentes del principio de equivalencia. Para ello, notemos que si se quiere ir del polo sur al polo norte de una esfera con desgaste mínimo, se puede recorrer cualquiera de sus meridianos. Por tanto, dos personas viajando simultáneamente por dos meridianos diferentes verán como su distancia mutua varía, creciendo inicialmente y disminuyendo después de cruzar el ecuador, hasta que se crucen de nuevo en el polo norte, y todo ello a pesar de que ambas siguen curvas geodésicas. Esta *aceleración relativa* que sufren cuerpos que se mueven *geodésicamente* es una consecuencia de la *curvatura* de la superficie de la esfera. De manera análoga, si dos cuerpos se dejan caer libremente desde una misma altura en el campo gravitatorio terrestre, es evidente que a medida que pasa el tiempo los dos cuerpos se acercan, dado que sus trayectorias apuntan hacia el centro de la Tierra. O sea, sufren una *aceleración relativa*. Está claro, por otro lado, que en un sistema de referencia no inercial, sin gravitación, dos partículas de prueba libres jamás sufren tal aceleración relativa. En resumidas cuentas, hay una manera de distinguir entre los campos gravitatorios *reales* y los sistemas de referencia no inerciales, *siempre y cuando se usen dos o más masas de prueba para sondear los efectos no locales.*

Pues bien, conjugando audazmente su principio de equivalencia, que llevaba a formas cuadráticas tales como (6), con las novedosas teorías matemáticas del genial Riemann sobre espacios curvos, Einstein alcanzó la maravillosa e inaudita conclusión de que la curvatura del espacio-tiempo era, pura y simplemente, el mismísimo campo gravitatorio. Dicho de otro modo, no hay fuerza de la gravedad, sino que simplemente las grandes masas deforman la geometría del espacio-tiempo y las partículas siguen las trayectorias extremales naturales de dicha geometría.

Para ello necesitó la inestimable ayuda de su amigo Grossmann. El concepto básico en la teoría tenía que ser el campo tensorial métrico fundamental $g$ que, aparte las restricciones mencionadas más arriba, era arbitrario y por lo tanto siempre se puede expresar como

$$g = \sum_{a,b=0}^{3} g_{ab}\, dx^a \otimes dx^b \equiv g_{ab}\, dx^a \otimes dx^b$$

en una cierta carta local $\{x^a\}$ de una variedad diferenciable cuadridimensional, donde $g_{ab} = g_{ba}$ son funciones sobre la variedad. Hay que notar que la segunda manera de escribir esa fórmula, sin el sumatorio, en la que se usa el convenio de sumación de Einstein, nació naturalmente en esta época de desarrollo de la nueva teoría. Un espacio con un tensor métrico fundamental tiene asociados, automáticamente, una única conexión canónica $\nabla$ –la conexión denominada de Levi-Civita, que satisface $\nabla g = 0$ y no tiene torsión–, y la correspondiente curvatura $\boldsymbol{R}$:

$$\left(\nabla_{\vec{X}}\nabla_{\vec{Y}} - \nabla_{\vec{Y}}\nabla_{\vec{X}} - \nabla_{[\vec{X},\vec{Y}]}\right)\vec{Z} \equiv \boldsymbol{R}(\vec{X},\vec{Y})\vec{Z}\,. \tag{6}$$

Aquélla define las trayectorias "geodésicas", que son a su vez extremales de la longitud $L$ definida como

$$L = \int_{\lambda_i}^{\lambda_f} \sqrt{|g(\vec{v},\vec{v})|}\, d\lambda \tag{7}$$

para cualquier arco continuo de curva $\gamma(\lambda)$ que una los puntos $p = \gamma(\lambda_i)$ y $q = \gamma(\lambda_f)$, sea diferenciable a trozos, y cuyo campo vectorial tangente en los trozos diferenciables es $\vec{v}$ –naturalmente, la integral se compone de la suma de los trozos diferenciables concatenados–. Estas geodésicas son las trayectorias que describen los cuerpos de prueba, o sea, los no sometidos a influencia externa alguna aparte de la gravedad. Y esto tanto si $g$ describe un campo gravitatorio real como si es la manifestación de fuerzas no inerciales. Obsérvese que hay curvas tales que $L = 0$ (las curvas luminosas, por las que se mueve la radiación electromagnética como la luz), y que por ello es evidente que $L$ no puede minimizarse para puntos con separación temporal[5]. En este caso $L$ puede, no obstante, maximizarse –lo que lleva directamente a la famosa "paradoja" de los gemelos, que no es tal, claro está–. Por su parte, el tensor de curvatura $\boldsymbol{R}$ mide, a través de la famosa ecuación de los campos de Jacobi, o igualmente de la llamada identidad de Ricci (6), la aceleración o desviación relativa de geodésicas cercanas. De esta manera, si $\boldsymbol{R} = 0$ no hay aceleración relativa, lo que sirve para describir sistemas no inerciales, mientras que si $\boldsymbol{R} \neq \boldsymbol{0}$ sí que la hay y se tiene un modelo de un cierto campo gravitatorio real. Como se sabe, matemáticamente la condición

$$\boldsymbol{R} = \boldsymbol{0}$$

identifica el caso en que el espacio-tiempo es *plano*, o sea, sin curvatura, y un teorema fundamental de Riemann asegura que se pueden escoger localmente coordenadas cartesianas tales que $g$ adopta la forma (4) si y sólo si $\boldsymbol{R} = 0$. El anterior edificio matemático queda redondo si además notamos, como hizo Einstein, que a lo largo de cualquier curva se pueden escoger coordenadas tales que los símbolos de Christoffel de la conexión se anulan, o sea, de tal forma que las geodésicas parecen ser, localmente, rectas. Esto es una versión matemática del principio de equivalencia local.

Del principio de equivalencia y la geometrización de la gravitación, que pasa a ser la curvatura del espacio-tiempo, se pueden deducir, directamente, algunos de los resultados más sorprendentes predichos por Einstein, tales como el cambio de la frecuencia de la luz al subir/bajar por la vertical (y en general el efecto Doppler gravitatorio) –véase (18) más adelante–, la desviación de la luz al pasar rasante al Sol, y en general la influencia del campo gravitatorio sobre la luz y sus trayectorias. Todo esto lo dedujo Einstein en los años 1907-13.

### 3.3 Las ecuaciones de campo de Einstein

Quedaba por resolver, empero, el problema más difícil. La fórmula (1) perderá su validez excepto como caso límite, ya que las fuerzas gravitatorias se han desvanecido y sólo queda la pura geometría. De forma similar, la ecuación de Poisson

$$\Delta\Phi = 4\pi G\varrho, \tag{8}$$

siendo $\Delta$ el operador laplaciano, que permite determinar el potencial gravitatorio[6] $\Phi$ dada la densidad de masa $\varrho$ de la distribución material que crea el campo, también queda en suspenso: ha de ser generalizada, a la vez que recuperada en un caso límite. Ahora bien, ¿cómo sustituir (1) y (8)? ¿Cómo cuantificar la curvatura del espacio-tiempo, la deformación producida por una distribución de masas? ¿Cómo determinar el tensor métrico $g$ en casos concretos de interés físico?

[5]Una cuestión que suele sorprender al profesional de las matemáticas no avisado es que una variedad lorentziana *no* es un espacio métrico en el sentido clásico, ya que no se puede definir una noción de distancia clásica. Por ello el clásico teorema de Hopf-Rinow no es válido, y de ahí que pueda haber muchos espacio-tiempos geodésicamente incompletos.

[6]El potencial $\Phi$ contiene la información de la energía por unidad de masa de un campo gravitatorio newtoniano, y la fuerza por unidad de masa correspondiente se obtiene mediante $\vec{f} = -grad\,\Phi$.

Einstein no supo resolver este problema durante años y, de hecho, propuso diversas ecuaciones para la nueva teoría. En un realmente ajetreado mes de Noviembre de 1915, en el que semanalmente presentó nuevas propuestas para las ecuaciones, finalmente se convenció de cuáles eran las correctas y las presentó el día 25. Se cree que este esfuerzo contribuyó a un deterioro gravísimo de su salud que por poco le hizo fenecer en 1917.

Para encontrar las *ecuaciones de Einstein*, que así es como se conocen desde entonces, el genio tuvo que relacionar ideas de índole diversa y usar una profunda intuición junto con el conocimiento que había adquirido de la geometría riemanniana. Dado que la teoría era relativista, la generalización de (8) debía contener también la constante fundamental $c$. Por otro lado, la relatividad especial había unificado la masa y la energía, y también ésta con la cantidad de movimiento en un 4-vector momento. Por ello, no sólo la densidad de masa sino todas las densidades de energía y momento debían aparecer como fuentes del campo gravitatorio. O sea, a la derecha de las ecuaciones. ¿Qué poner a la izquierda? Para descifrar este misterio se podía usar como inspiración el campo métrico válido en la aproximación newtoniana, que Einstein ya creía conocer basándose en su principio de equivalencia. Este campo $g$ se puede escribir, en un sistema inercial tradicional, como

$$g = -(c^2 - 2\Phi)\, dt \otimes dt + dx \otimes dx + dy \otimes dy + dz \otimes dz \qquad (9)$$

por lo que se intuye que las funciones $g_{ab}$ serán, en general, las que sustituyan al potencial $\Phi$ (por ello se denomina a veces potenciales a las $g_{ab}$). Obsérvese que la generalización no es baladí . . . , se pasa ¡de una a diez funciones! Se sigue de ello que las derivadas primeras (los símbolos de Christoffel de la conexión $\nabla$) son análogos a las fuerzas gravitatorias de Newton, lo que concuerda con la expresión en coordenadas de la ecuación de las geodésicas. Y teniendo en cuenta que la relación (8) es lineal y contiene las derivadas segundas de $\Phi$, se espera por lo tanto que las buscadas ecuaciones tengan en el miembro izquierdo *algo* que contenga las derivadas segundas de las $g_{ab}$, a ser posible que sea lineal en ellas.

¿Qué puede ser ese algo? Para empezar, es necesario que esté matemáticamente bien definido, o sea, que sea un objeto con carácter tensorial para poder considerarlo sin problemas en bases arbitrarias. Una posibilidad sería poner un campo vectorial proveniente de las derivadas segundas de $g$, para igualarlo a la densidad de 4-momento de la distribución material y energética. No obstante, no era posible construir un campo vectorial de las citadas características. El único objeto tensorial que se puede fabricar con las derivadas segundas de $g$, lineal en éstas, es ni más ni menos que $\boldsymbol{R}$. Pero el tensor de curvatura es demasiado complicado, y tiene un excesivo número de componentes independientes (20 en dimensión 4). A lo mejor que se podía aspirar, por lo tanto, era a usar una traza de $\boldsymbol{R}$. La traza natural es el conocido tensor de Ricci $\boldsymbol{Ric} \equiv \mathrm{tr}\boldsymbol{R}$, o usando notación de índices

$$R_{ab} \equiv R^{c}{}_{acb}\,,$$

donde $R^{a}{}_{bcd}$ representa el tensor de Riemann asociado a la curvatura $\boldsymbol{R}$. $\boldsymbol{Ric}$ es un campo tensorial 2-covariante simétrico. A primera vista se puede pensar en igualar $\boldsymbol{Ric}$ al tensor energía-momento de la materia $\boldsymbol{T}$, que es del mismo tipo y ya se conocía bien a partir de los trabajos sobre el campo electromagnético y también de los estudios en dinámica de fluidos. No obstante, la relación $\boldsymbol{Ric} \propto \boldsymbol{T}$ tenía graves problemas, uno de ellos es que no se deduce el límite newtoniano correcto, o sea, no conduce a (8) usando (9). Otro grave problema es que estas ecuaciones no incorporan las necesarias relaciones de conservación de la energía y el momento de la materia, que se pueden expresar como

$$\boldsymbol{\nabla} \cdot \boldsymbol{T} = 0, \tag{10}$$

donde $\boldsymbol{\nabla}\cdot$ denota la divergencia. Hay que subrayar, para ser fieles a la verdad, que Einstein no usó esta propiedad y no la mencionó en sus trabajos hasta después de haber encontrado las ecuaciones correctas.

Finalmente, las ecuaciones que propuso, y que excepto por lo que se dirá en la próxima sección siguen hoy plenamente vigentes, fueron

$$\boldsymbol{Ric} - \frac{\mathrm{tr}\boldsymbol{Ric}}{2}\, g = \frac{8\pi G}{c^4}\boldsymbol{T}, \tag{11}$$

donde $\mathrm{tr}\boldsymbol{Ric}$ denota la traza respecto de $g$ de $\boldsymbol{Ric}$, o sea, la *curvatura escalar*. El miembro izquierdo de (11) es un campo tensorial 2-covariante simétrico que se denomina *tensor de Einstein* y en virtud de las identidades de Bianchi para la curvatura satisface *idénticamente* la propiedad

$$\boldsymbol{\nabla}\cdot\left(\boldsymbol{Ric} - \frac{\mathrm{tr}\boldsymbol{Ric}}{2}\, g\right) = 0,$$

recuperando así de manera manifiesta la ley (10) en toda situación posible.

Como se aprecia, las ecuaciones de Einstein relacionan las propiedades métricas del espacio-tiempo junto con su deformación en términos de la curvatura con la cantidad y distribución energético-material que crea esa deformación. La fórmula (11) representa un sistema de diez ecuaciones diferenciales en derivadas parciales no lineales y acopladas de segundo orden, por lo que son matemáticamente muy complicadas. Desde un punto de vista físico no se conocía nada igual, ya que las fuentes crean un campo gravitatorio que es una deformación de la geometría subyacente, pero la propia materia se tiene que mover según dicta dicha geometría, y este movimiento afecta a la geometría de nuevo que a su vez modifica la dinámica de las fuentes..., y así sucesivamente en un baile sin fin de la pescadilla que quiere morderse la cola.

## 4 1917: el nacimiento de la Cosmología teórica

El año 1916 se publicó el artículo fundacional de la Relatividad General. Se había dado así un paso de gigante en la física, que fue además fundamental para el nacimiento de la Cosmología en una doble vertiente: por disponer de una nueva teoría de la gravedad, y porque ésta relacionaba la materia con la geometría del espacio y del tiempo. Tales herramientas en manos del genio de Einstein eran las piezas de una "bomba de relojería", valga la expresión, por lo que la eclosión de una ciencia del cosmos era inevitable. Además, la ciencia contó en este caso con un gran catalizador e instigador, Wilhem de Sitter, astrónomo holandés que discutió con Einstein en 1916 acerca de la nueva teoría y de la posibilidad de que la inercia tuviera un origen totalmente material –una idea machiana como ya se ha mencionado–, lo que parecía requerir un Universo *finito*. Por cierto que de paso esto parecía arreglar la paradoja de Seeliger.

Estas discusiones influyeron en Einstein decisivamente, de manera que en 1917 publicó el articulo fundacional de la Cosmologia teórica, titulado *Kosmologische Betrachtungen zur allgemeinen Relativitätstheorie* [1]. En esta contribución fundamental se presentaba un modelo de Universo conocido hoy día como el *universo estático de Einstein*.

En la época que nos ocupa, la mayoría pensaba que el Universo era *estático* y consistía en la Vía Láctea y si acaso vacío más allá. Esta creencia es fácil de entender, de aceptar como sensata y basada en las observaciones, si pensamos por un lado que –a pesar de la premonitoria idea kantiana de los universos isla– las galaxias como tales no se habían descubierto todavía, y por otro que los movimientos estelares *visibles* son claramente periódicos y debidos a la rotación de la Tierra sobre sí misma y en torno al Sol: nadie había pensado en movimientos radiales, de alejamiento, ni era lógico ni fácil buscarlos. Einstein quiso *modelizar el Universo entero* y para ello supuso que, en buena aproximación, todas las estrellas tenían velocidades despreciables respecto de un sistema de referencia adecuado de manera que la distribución material se podía considerar, a parte de estática, *espacialmente homogénea e isótropa.* Las razones para esta hipótesis eran parcialmente observacionales, pero tampoco parece muy arriesgado aventurar que las ideas de simplicidad, sencillez, democracia y erradicación de privilegios inspiraron al genio una vez más. Con ello dio lugar a un nuevo principio de "humildad": no ocupamos un lugar especial en el Universo. Y como en tantas otras ocasiones, este paso de modestia desencadenó un salto adelante gigantesco. La estaticidad y este principio le condujeron a un campo métrico del tipo

$$g = -c^2 dt \otimes dt + {}^3g$$

donde ${}^3g$ es el tensor métrico tridimensional de cualquier espacio homogéneo e isótropo, o en matemáticas palabras, un espacio *máximamente simétrico.*

Aludiendo aquí a las ideas machianas que tanto le inspiraron, y tratando de evitar la desagradable tarea de decidir cuales debían ser las condiciones de contorno adecuadas para *todo el Universo*, Einstein consideró que el espacio debía ser *finito pero ilimitado.* Lo cual conduce a que ${}^3g$ sea el tensor métrico natural de una 3-esfera (o sea, el ${}^3g$ heredado por una esfera tridimensional inmersa de manera natural en $\mathbb{R}^4$). Naturalmente, una tal 3-esfera tiene un volumen finito, no hay límites, y todos sus puntos son equivalentes. Y esta simple idea fue otro golpe de intuición sin igual, una drástica revolución en la física y un éxito de las matemáticas aplicadas. Es el primer caso de la historia –de la física, o sea del mundo *real*– en que *se cambia la topología del espacio.* En lenguaje matemático esto equivale a decir que se ha cambiado la variedad subyacente. Así, el espacio-tiempo es la variedad $\mathbb{R} \times S^3$ con campo métrico

$$g = -c^2 dt \otimes dt + a^2 \left[ d\chi \otimes d\chi + \sin^2 \chi \left( d\theta \otimes d\theta + \sin^2 \theta \, d\phi \otimes d\phi \right) \right] \qquad (12)$$

en una carta local clásica donde $a$ es una constante que representa el radio de curvatura de la 3-esfera de forma que su volumen total es $2\pi^2 a^3$.

Surgió entonces un problema inesperado y aparentemente irresoluble que puso a prueba el inigualable coraje de Einstein. Resulta que las ecuaciones de campo (11) no admiten *ninguna* solución estática a la vez que espacialmente homogénea e isotrópa. Einstein no se arredró, desechó la idea de abandonar la isotropía y homogeneidad del espacio universal a gran escala (¡qué sabiduría!), y ni siquiera tomó en consideración la posibilidad de renunciar a un mundo estático[7]. Desapasionadamente, Einstein decidió que no había problema, el fin justifica el camino, y pensó: "simplemente", ¡cambiemos las ecuaciones! Después de años detrás de las ecuaciones (11), y cuando ya todo parecía cuadrar, decidió desdecirse, y ello aunque hubiera de corregirse a sí mismo en la meta que más le había costado conseguir y que casi le privó de una buena salud. Y así fue cómo las ecuaciones (11), que tantos esfuerzos le ocasionaron, fueron modificadas con tal de acomodar una solución cosmológica estática y uniforme. Las ecuaciones modificadas toman la forma

$$\boldsymbol{Ric} - \frac{\mathrm{tr}\boldsymbol{Ric}}{2} g + \Lambda g = \frac{8\pi G}{c^4} \boldsymbol{T} \qquad (13)$$

donde el nuevo sumando $\Lambda \boldsymbol{g}$ se denomina "término cosmológico" y la nueva constante $\Lambda$ se llama "constante cosmológica". Obviamente, la condición de nula divergencia para el término de la izquierda sigue siendo válida

$$\boldsymbol{\nabla} \cdot \left( \boldsymbol{Ric} - \frac{\mathrm{tr}\boldsymbol{Ric}}{2} g + \Lambda g \right) = 0$$

[7]Hoy podemos decir que esto fue una obstinación, una clara tozudez... Pero ¿quién en su sano juicio hubiera dicho en 1917 que el Universo no era estático? En este sentido, véase más adelante la discusión de la magnífica aportación de Friedman y la reacción de Einstein.

por lo que la ecuación fundamental de conservación (10) se mantiene plenamente. Ahora se sabe que, de hecho, el lado izquierdo de la relación (13) es el más general posible compatible con esta conservación junto con su tensorialidad y la linealidad en las derivadas segundas de $g$.

Para apreciar mejor la radicalidad de la introducción de $\Lambda$, notemos que con el nuevo término cosmológico la atracción gravitatoria en el límite de bajas velocidades, o sea en el límite newtoniano, se puede describir aproximadamente con la siguente modificación de la fórmula (1)

$$\vec{F} = \left( \frac{GMm}{r^2} - \frac{\Lambda}{3} mc^2\, r \right) \frac{\vec{r}}{r} \tag{14}$$

donde constatamos que, para distancias cortas, se recupera la ley de Newton; pero para largas distancias domina el segundo sumando que es repulsivo (si $\Lambda$ es positiva). Esto de hecho impone estrictos límites experimentales en los valores de una $\Lambda$ que sea constante universal si se quiere recuperar la astronomía planetaria del sistema solar.

En el artículo en consideración se presentó una solución de las nuevas ecuaciones (13) verificando las hipótesis adoptadas. Curiosamente, la resolución de las ecuaciones (13) para el campo (12) cerraba el círculo de manera elegante, redondeaba la tarea cuantificando la relación entre la cantidad de materia y la forma del Universo mediante la fórmula exacta

$$\Lambda = \frac{1}{a^2} = \frac{4\pi G}{c^2}\, \varrho \tag{15}$$

siendo $\varrho$ la densidad de masa del cosmos –que es constante en este el primer modelo del mismo– para una materia enrarecida que no siente presiones ni tensiones ni contiene flujos de energías:

$$\boldsymbol{T} = \varrho c^4 dt \otimes dt\,. \tag{16}$$

La expresión (15) implica que para universos pequeños la densidad ha de ser grande y, viceversa, si la 3-esfera es muy grande la densidad de masa es muy reducida.

Finalmente, Einstein analizó la cuestión de si la inercia de las partículas de prueba puede ser una manifestación de la influencia de todo el Universo, y llegó a la conclusión de que, en este modelo, la inercia está influida, pero no determinada completamente, por la materia universal. Estas conclusiones son dudosas y actualmente es inusual perseguir las ideas machianas en el marco de la relatividad general. De hecho, Einstein fue perdiendo poco a poco interés en estas ideas hasta rechazarlas totalmente en 1954 [2].

Resumiendo lo que se ha escrito en esta sección, podemos recopilar en la tabla siguiente las ideas germinales, atrevidas y geniales que Einstein presentó de forma natural en su artículo, erigiendo una instigadora lista que irremediablemente cautiva.

| |
|---|
| **El Universo como objeto de la física.**<br><br>Por tanto, nacimiento de la ciencia de la Cosmología teórica. |
| **Homogeneidad e isotropía espacial.**<br><br>Esta idea se conoce ahora como el *Principio Cosmológico.*<br>En boga actualmente y con amplio apoyo observacional.<br>Es la base de los modelos cosmológicos corrientes. |
| **Estudio de la 'forma' del Universo.**<br><br>Cambio de la topología del espacio por primera vez en la historia.<br>Esto es un ejercicio habitual en la ciencia actual. |
| **Propuesta de un espacio finito e ilimitado.**<br><br>Esta posibilidad sigue abierta.<br>Parece ser cierta si el Universo es 'cerrado'—ver sección 6. |
| **Relación precisa entre la geometría del Universo y su distribución material.**<br><br>Naturalmente, esto es una consecuencia directa de la relatividad general.<br>Por ello, continúa (y continuará) vigente adecuadamente modificada. |
| **Corrección de sus ecuaciones originales; introducción del término cosmológico.**<br><br>Ahora hay muchas razones para creer que $\Lambda$ no es nula.<br>Gran controversia actual. |
| **Búsqueda de una explicación científica de la inercia.**<br><br>Obsoleta. |

# 5 Contribución de los matemáticos entre 1917 y 1933

Como ya se ha mencionado, de Sitter influyó decisivamente en el despegue de la ciencia cosmológica por medio de sus conversaciones con Einstein. De mayor importancia probablemente fueron, si cabe, sus contribuciones directas al nacimiento y desarrollo de las ideas cosmológicas en la relatividad, mediante dos artículos fundamentales publicados el mismo 1917 titulados: *On the Relativity of Inertia, remarks concerning Einstein's latest Hypothesis* y *On the Curvature of Space* [3].

## 5.1 La contribución de de Sitter

Willem de Sitter se había licenciado en matemáticas en Groningen, por lo que tuvo una buena formación en esta rama de la ciencia a pesar de que se doctoró y destacó como astrónomo, siéndole concedida la medalla Bruce en 1931. Su sapiencia matemática le permitió comprender y usar sin mayores problemas las herramientas de geometría subyacentes a la relatividad general, de manera que en los artículos mencionados presentó y discutió una solución de las ecuaciones modificadas (13) que representa un modelo totalmente *vacío*, o sea, tal que el miembro derecho de (13) se anula, lo que significa ni más ni menos que no hay ningún tipo de materia[8]. El campo tensorial métrico de la solución de de Sitter es

$$g = -c^2 \cos^2 \chi \, dt \otimes dt + a^2 \left[ d\chi \otimes d\chi + \sin^2 \chi \left( d\theta \otimes d\theta + \sin^2 \theta \, d\phi \otimes d\phi \right) \right] \quad (17)$$

por lo que parece compartir la variedad subyacente, y la 3-esfera espacial, con (12). No obstante, hay que resaltar que ahora el rango permitido de valores de $\chi$ no es el habitual (evitando las singularidades en los polos $0 < \chi < \pi$), porque para $\chi \to \pi/2$ aparece una singularidad en $g$, que deviene un tensor métrico degenerado con radical no vacío. Esta propiedad indica la existencia de un *horizonte*, un concepto de la máxima importancia en relatividad y que desafortunadamente no podemos analizar aquí ni medianamente. Baste decir, para los lectores con aficiones físicas, que el tipo de comportamiento que aparece en $\chi \to \pi/2$ es análogo al visto en (6) para $\omega\rho \to c$, el cual es intuitivamente más asequible. En realidad, es muy sencillo comprobar que (17) es (parte de) un espacio-tiempo –es decir como variedad cuadridimensional– máximamente simétrico y por ello de curvatura constante (positiva). Es fácil entonces encontrar cambios de coordenadas que *extienden* (17) a todo el espacio-tiempo global con curvatura constante y positiva –véase ( 5.4 ) más abajo–, lo que es análogo al cambio que extiende (6) a todo el espacio-tiempo plano de Minkowski (4).

[8]La solución de de Sitter es muy importante en aplicaciones y estudios actuales, debido a que presenta una etapa inflacionaria y un horizonte –véase el apartado 5.4 más abajo–, conceptos que eran desconocidos o malentendidos en 1917, pero que contribuyó a aclarar y explicar.

En todo caso, calculando la curvatura de (17) se puede llegar rápidamente a la conclusión de que este campo métrico es una solución de (13) con

$$\boldsymbol{T} = \mathbf{0}, \qquad \Lambda = \frac{3}{a^2}$$

de donde se deduce que *no hay materia* en este segundo universo estático. Esta indeseada propiedad cambió radicalmente el comienzo incipiente de la Cosmología, ya que prueba que las partículas de prueba, que se mueven por las geodésicas temporales de (17), experimentan una inercia que, habida cuenta de que el Universo está vacío, no puede ser machiana, sino que parece ser *relativa al espacio*. Un resultado que hizo desmoronarse el esquema einsteinano sobre el cosmos, la inercia, el principio de Mach y les ecuaciones modificadas (13). Seguramente en este preciso momento empezó a arrepentirse de haber cambiado sus ecuaciones originales (11)...

De mayor relevancia aún para lo que nos ocupa es otra inesperada propiedad de la solución (17): imaginemos una partícula cualquiera en reposo en un punto arbitrario del espacio $(\bar{\chi}, \bar{\theta}, \bar{\phi})$. Si calculamos $L$ para esta curva del espacio-tiempo entre dos valores de $t$, según la fórmula (7) obtenemos

$$\frac{L_{12}}{c} = (t_2 - t_1) \cos \bar{\chi}$$

lo que da el tiempo propio transcurrido para estas partículas entre $t_1$ y $t_2$. Como se ve, este tiempo propio *depende* de la posición de la partícula. Por ello, dos personas distintas viviendo en este universo ven que sus relojes marchan a un ritmo diferente según la posición que ocupen en (el trozo de) la 3-esfera. Naturalmente, esto se puede medir, por ejemplo emitiendo uno onda electromagnética (un rayo de luz) con frecuencia fija desde un cierto punto y recibiéndola en otro. Al cambiar la medida del tiempo con $\chi$, la onda recibida tiene una *frecuencia distinta* de la de la onda emitida. Esto se mide mediante el parámetro $z$ de desplazamiento ("shift" en inglés)

$$1 + z \equiv \frac{\nu_e}{\nu_r} = \frac{\lambda_r}{\lambda_e}$$

siendo $\nu$ la frecuencia y $\lambda$ la correspondiente longitud de onda para el emisor ($e$) y el receptor ($r$). No hay que olvidar que este tipo de efecto había sido calculado ya por Einstein mucho antes usando la aproximación newtoniana (9), de donde se sigue fácilmente

$$1 + z = \sqrt{\frac{c^2 - \Phi_e}{c^2 - \Phi_r}} \qquad \Longrightarrow \qquad z \approx \frac{1}{2c^2}(\Phi_e - \Phi_r)\,, \tag{18}$$

pero nunca en un contexto cosmológico. Incluso de Sitter rebuscó entre los resultados experimentales para saber si había fuentes estelares lejanas (insistamos en que las galaxias no se conocían aún) de las que recibiéramos sus espectros desplazados hacia el rojo o el azul. Desafortunadamente, las medidas de la época al respecto eran escasísimas y los primeros resultados concluyentes de Slipher (véase la sección 5.3 ) fueron publicados en 1915-17.

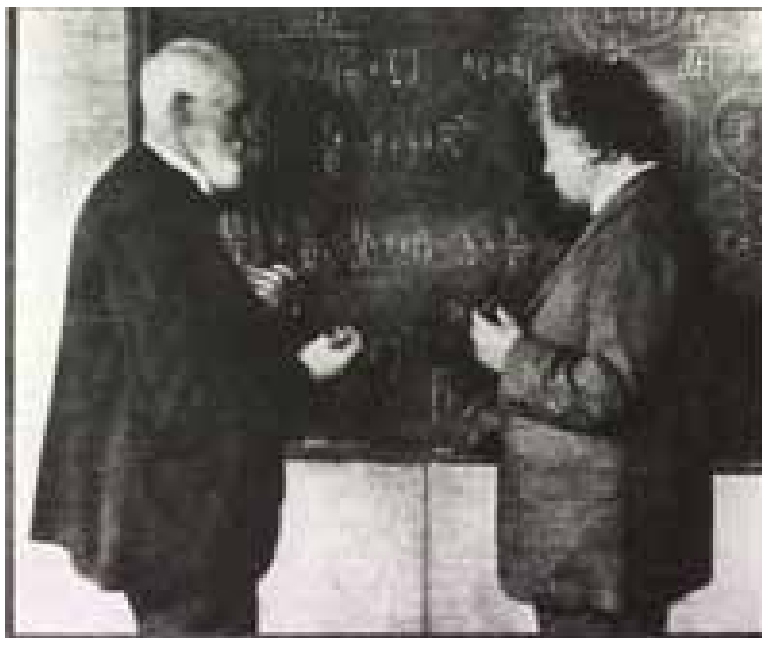

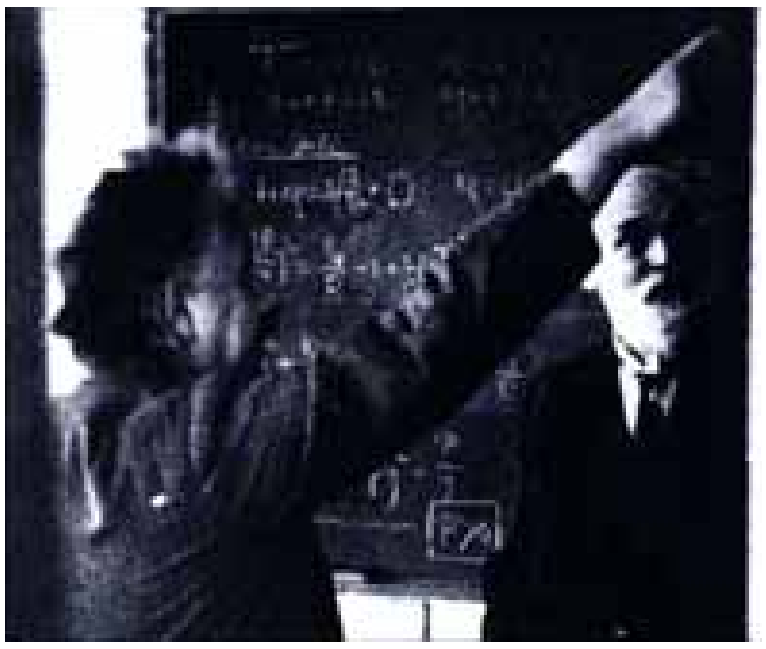

A. Einstein y W. de Sitter discutiendo en el *International School of Advanced Science* del *California Institute of Technology*, Pasadena, en enero de 1932.

## 5.2 La osadía matemática de Friedman

El siguiente paso esencial, que además portaría al poco tiempo a la síntesis definitiva, fue llevado a cabo en 1922 por el matemático ruso Aleksandr Aleksandrovich Friedman (también se usa a menudo Friedmann, germanización de la transcripción directa del ruso cirílico original 'Fridman'), quien en un artículo de relevancia histórica titulado *Über der Krümmung des Raumes* [4] presentó una nueva familia de soluciones de las ecuaciones (13) de acuerdo con el principio cosmológico, es decir, para distribuciones materiales espacialmente homogéneas e isótropas. La historia de este artículo es ejemplar desde el punto de vista de la relación entre las matemáticas y la física, y tiene una moraleja de hondo calado para los profesionales de aquélla: la segunda avanza muchas veces gracias a las aportaciones originalísimas de los matemáticos.

El planteamiento de Friedman no puede ser más lógico y natural para un matemático. Se trata de encontrar la familia de soluciones de las ecuaciones (13) que cumplen ciertas hipótesis (un espacio de curvatura constante positiva, materia sin presión ni flujos) deduciendo de paso las soluciones conocidas de Einstein y de Sitter e identificándolas como las que satisfagan alguna propiedad dentro de la familia general, si es que ésta existe (en otro caso, se probaría la unicidad de las dos soluciones mencionadas). Por otro lado, el planteamiento de Friedman no podía sonar más descabellado a cualquier físico razonable de la época, ya que permitió que la curvatura constante del espacio *dependiera* *del tiempo*, en flagrante contradicción con el *dogma* vigente: el Universo era estático.

En fin, el caso es que Friedman probó que con esas hipótesis el tensor métrico debía adoptar necesariamente la forma

$$g = -M^2\, dt \otimes dt + a^2(t) \left[ d\chi \otimes d\chi + \sin^2 \chi \left( d\theta \otimes d\theta + \sin^2 \theta\, d\phi \otimes d\phi \right) \right] \quad (19)$$

donde $M$ es una función cualquiera en $\mathbb{R} \times S^3$ y ahora $a(t)$ depende de $t$. Y demostró entonces que las ecuaciones (13) se satisfacen para un tensor $\boldsymbol{T}$ del tipo (16) solamente si $\dot{a}\, dM = 0$, donde un punto indica derivada respecto de $ct$, lo que lleva a dos posibles situaciones.

1. $a$ es constante. Entonces el resto de las ecuaciones conducen irremediablemente a una de las dos soluciones conocidas (12) o (17). De este modo las soluciones conocidas quedan definitivamente caracterizadas como las que son estáticas dentro de la familia general.

2. $\dot{a} \neq 0$. En este caso se sigue que $M$ ha de ser constante y siempre se puede escoger $M = c$.

En este segundo caso el conjunto completo de ecuaciones de Einstein conducen a

$$\frac{8\pi G}{c^2}\varrho + \Lambda = \frac{3}{a^2}(\dot{a}^2 + 1), \qquad \Lambda = 2\frac{\ddot{a}}{a} + \frac{1}{a^2}(\dot{a}^2 + 1)$$

lo que proporciona una relación precisa entre *la evolución* de las propiedades geométricas del espacio-tiempo y las variables de la materia. La segunda de estas ecuaciones tiene una integral primera inmediata

$$a(\dot{a}^2 + 1) = A + \frac{\Lambda}{3}a^3$$

para una constante $A$, lo que permite resolver para $a(t)$ siquiera implícitamente en términos de integrales elípticas. La otra ecuación da entonces la forma explícita de la densidad de masa

$$\varrho = \frac{3c^2}{8\pi G}\,\frac{A}{a^3}\,.$$

Nótese que la constante cosmológica queda libre, como un parámetro, y que por lo tanto es innecesaria: existe una solución de Friedman para las ecuaciones originales (11).

Las soluciones de Friedman estaban en *expansión o contracción* dependiendo del signo de $\dot{a}$. En todo caso eran dinámicas, lo que conllevaba un parámetro de desplazamiento $z$ no nulo. Pero el análisis de las soluciones para $a(t)$ llevaba a la conclusión de que en una mayoría de casos (para $\Lambda < 4\pi G\varrho/c^2$) era inevitable que $a \to 0$ para un valor de $t$ finito, en el "futuro" si $\dot{a} < 0$,

o en el pasado en el otro caso. Matemáticamente esto es una singularidad catastrófica, ya que el campo tensorial métrico pasa a tener un radical tridimensional (el espacio) y la curvatura $\boldsymbol{R}$ diverge. Físicamente esto es peor que una catástrofe, ¡es el fin del mundo!— casi literalmente. El espacio desaparece y la densidad de masa se hace infinita. Friedman hablaba en su artículo del "tiempo de creación" para designar el valor de $t$ desde $a \to 0$ hasta el instante actual si $\dot{a} > 0$.

Si aquello de permitir que el Universo fuera dinámico ya era malo, la catástrofe anterior era el colmo, no sólo inaudita, sino que a cualquier científico de la época le debió parecer una extravagancia inaceptable. La reacción de Einstein fue inmediata, fulminante y, claro está, negativa. De hecho, las hipótesis y las consecuencias de los modelos de Friedman eran tan descabelladas que Einstein publicó casi instantáneamente una réplica [5] asegurando que tales modelos no eran soluciones de (13): "Los resultados acerca de un mundo no estacionario [...] me parecen sospechosos. En realidad, resulta que la solución dada no satisface las ecuaciones de campo" [5]. Claro está, en este caso erró de plano. Al cabo de poco, Friedman logró convencerle, epistolarmente, de que sus cálculos eran correctos. La historia tiene su miga, porque Friedman escribió a Einstein el 6 de Diciembre de 1922:

> "Considerando que la posible existencia de un mundo no estacionario tenga cierto interés, me permito presentarle aquí los cálculos que he hecho [...] para su verificación y evaluación crítica. [...] Si encontrase estos cálculos correctos, le ruego por favor que sea tan amable de informar de ello a los editores de *Zeitschrift für Physik* [...]".

La carta llegó a Berlin cuando Einstein ya no estaba allá, y adonde no volvería hasta el 15 de Marzo de 1923. ¿Dónde estaba? De gira por Japón, Palestina, y finalmente degustando caldos y viandas sin igual en... ¡España! [6]. Ni siquiera a su vuelta se percató o tuvo tiempo de leer la misiva de Friedman, de la que parece tuvo conocimiento cuando un colega de éste llamado Krutkov alertó a Einstein, en Mayo de ese año, acerca de ella y su contenido. Sólo entonces admitió Einstein su error y rápidamente lo corrigió en una nueva nota [7, 8]:

> "En mi nota previa criticaba [4]. Sin embargo, mi crítica, como me ha convencido la carta de Friedman que me trasladó el Sr. Krutkov, se basaba en un error de mis cálculos. Creo que los resultados del Sr. Friedman son correctos e iluminadores ... "

de la que sólo en el último momento (¿buena educación?, ¿intuición?) eliminó una frase donde literalmente decía que los resultados de Friedman eran *matemáticamente correctos, pero físicamente irrelevantes* [8]. Einstein simplemente no podía imaginarse la posibilidad de que el Universo no fuera estacionario. Ahora podemos decir que él y el resto de científicos estaban obsesionados –eso sí, basándose en observaciones– con la estacionariedad del cosmos. No había pruebas, ni siquiera indicios, de que el Universo sufriera una evolución dinámica. La predicción teórica de la expansión universal quedó así, desafortunadamente, frustrada de nuevo.

En resumidas cuentas, Friedman recogió todas las ideas geniales de Einstein (y de Sitter) pero descartando, o al menos abriendo la posibilidad de descartar, la inmutabilidad del cosmos. Me gustaría resaltar solemnemente que los matemáticos tienen en ocasiones, muchas más veces de lo que podamos pensar, esta gran ventaja: su falta de "prejuicios físicos" les permiten descrubrir resultados que están vedados a los mejores físicos. Esta es la moraleja. Sospecho, por otro lado, que la falta de "prejuicios matemáticos" de estos últimos ha permitido devolver el favor con creces.

Aunque pueda parecer que los prejuicios sobre la estacionariedad del cosmos eran razonables –que lo eran– y que por lo tanto la historia siguió un curso hasta cierto punto lógico, hay que señalar, mal que le pese al autor de estas líneas, que en este caso concreto los físicos teóricos no estuvieron a la altura de las circunstancias. Digo esto porque cualquier físico medianamente formado, no digamos ya el genial Einstein, se debería haber percatado inmediatamente, al verse confrontado con las soluciones de Friedman, de que el universo estático de Einstein es *inestable*: cualquier mínima perturbación conduciría a una situación dinámica, a su colapso irremediable. De hecho, este fue uno de los argumentos posteriores de Lemaître y Eddington –véase más abajo–. Pero ya se sabe, así se escribe la historia.

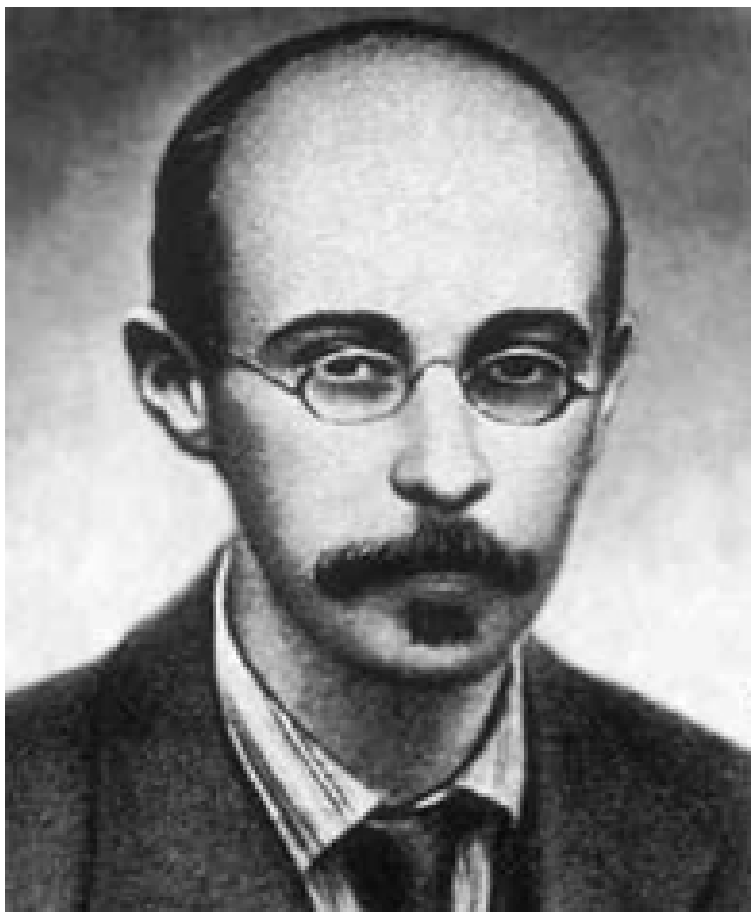

A.A. Fridman

## 5.3 Un paréntesis astronómico: las observaciones

A pesar de que puede ser de menor interés para un matemático, necesitamos ahora un breve paréntesis sobre la revolución que tuvo lugar, en paralelo, acerca de las dimensiones y constitución del Universo en el periodo 1915-30. Esto es así porque fueron las observaciones astronómicas las que destruyeron final y definitivamente el prejuicio de la inmutabilidad del Universo. El final feliz de esta historia llegó cuando los desarrollos teóricos que se han ido relatando se vieron complementados por la experimentación, que como se sabe siempre tiene la última palabra en la ciencia.

La estructura del Universo, la distribución de la materia, y la distancia de los objetos cósmicos eran cuestiones de grandísima controversia en esa época. El momento álgido ocurrió el 26 de abril de 1920 con el renombrado "Gran debate de la Astronomía", que tuvo como principales protagonistas a Heber

Doust Curtis y Harlow Shapley. Se discutieron muchos temas en ese debate que se tituló *The Scale of the Universe*, nos interesa ahora su parte más conocida referente a los divergentes puntos de vista sobre la existencia de nebulosas espirales *extra-galácticas* –¡los universos-isla de Kant!–. Se pueden consultar fuentes más precisas en [9], pero en pocas palabras Curtis estaba a favor de su existencia, mientras que Shapley defendía lo contrario[9].

La resolución del debate se produjo con las aportaciones trascendentales del gran Edwin Powell Hubble en uno de los momentos clave de la historia de la ciencia. Usando el histórico telescopio Hooker del observatorio del Monte Wilson (California, EE UU) en los años 1923-24, Hubble descubrió estrellas cefeidas en los brazos espirales de la "nebulosa" M31 (Andrómeda), demostrando con ello que la distancia a M31 era mucho mayor que las estimaciones más optimistas del diámetro de la Vía Láctea. Por consiguiente, M31 estaba *fuera* de ésta, y era de hecho *otra* galaxia como la propia Vía Láctea. Esto mismo lo corroboró con las nebulosas M33 y NGC 6822, todas ellas miembros de nuestro *grupo local* de galaxias. Hubble estimó la distancia a M31 en 900000 años-luz (la distancia aceptada hoy en día es de aproximadamente 2,9 millones de años-luz), lo cual engrandecía el Universo conocido de manera insospechada, radical y explosiva.

Hubble también merece el crédito de haber medido el sistemático desplazamiento al rojo de los espectros recibidos de las nebulosas (hoy galaxias), y de la ley que lleva su nombre. El primer astrónomo que midió tales comportamientos fue, en realidad, Vesto Melvin Slipher, quien durante los años 1913-25 analizó los espectros de muchas nebulosas descubriendo una clara preferencia por el desplazamiento hacia el rojo: 11 nebulosas de 15 en 1915 [10], y una razón de 21 a 4 (rojo versus azul) en 1917 [11]. En este artículo Slipher defendió que, teniendo en cuenta que nosotros teníamos un movimiento respecto de las nebulosas[10], pero no respecto de las estrellas, las observaciones favorecían la idea de que las nebulosas espirales eran sistemas estelares a grandes distancias. Esto años antes de que Hubble descubriera las cefeidas en Andrómeda. Dejando esto a un lado, Hubble sí que descubrió la relación existente entre el desplazamiento hacia el rojo y la distancia a las nebulosas/galaxias. Esta relación es aproximadamente *lineal*

$$z \approx \frac{H_0}{c} D \qquad \text{(ley de Hubble)} \tag{20}$$

[9] Para hacer justicia, hay que decir que muchas de las otras proclamas de Curtis resultaron ser erróneas, y que por su parte Shapley fue el responsable de la tarea "copernicana" de poner al Sol (y de paso a la Tierra) en su sitio . . . , o sea, en los suburbios de la Vía Láctea –contra lo que Curtis pensaba–, y quien dio una estimación precisa de sus dimensiones.

[10] El desplazamiento al rojo se puede interpretar como efecto Doppler, o sea, como efecto de un movimiento de alejamiento o recesión de los objetos, análogamente a lo que ocurre con los pitidos de un tren cuando se aleja de la estación.

siendo $D$ la distancia al emisor, de manera que cuanto más lejana esté una galaxia, mayor es su velocidad de recesión. La constante $H_0$ se denomina *constante de Hubble* y tiene un valor controvertido que en el momento de escribir estas líneas parece ser $H_0 \simeq 70 \pm 5$ Km/s Mpc $= (13 \pm 1) \times 10^9$ años$^{-1}$. Esta es la base experimental de la expansión del Universo, por ello se considera uno de los cimientos de la Cosmología moderna, y es seguramente uno de los descubrimientos más inesperados y sorprendentes de la historia de la ciencia.

## 5.4 De nuevo los matemáticos: Lanczos y Weyl

Dos matemáticos ilustrísimos hicieron aportaciones de relevancia al tema de la expansión del Universo, el desplazamiento al rojo, y los modelos de Friedman y de Sitter. Primero, fue Cornelius Lanczos[11] quien en 1922 [12] escribió la solución de de Sitter –con más precisión habría que decir *una extensión métrica* de esa solución– en forma de modelo de Friedman. Seguidamente, dedujo la expresión para el desplazamiento al rojo en este modelo *dinámico* [13].

Algo parecido hizo independientemente Hermann Weyl en otro magnífico artículo [14]. Quizás ésta sea la primera discusión en profundidad de las fórmulas para el desplazamiento al rojo en modelos cosmológicos relativistas, y la fórmula obtenida se reduce, para velocidades pequeñas, a la expresión (20) que posteriormente vendría en llamarse ley de Hubble.

Hay que señalar además que, entre otras muchas contribuciones a la geometría riemanniana y la relatividad –como por ejemplo el tensor de curvatura conforme–, Weyl dio la primera versión rigurosa del Principio Cosmológico pero desde el punto de vista de un conjunto de observadores, o sea, de una congruencia diferenciable de curvas temporales en la variedad. O dicho sucintamente, usando un campo vectorial temporal. La inspiración vino al tratar de reconciliar una teoría covariante general (tensorial) como la relatividad general con la posibilidad de describir un único y particular conjunto de fenómenos físicos, el Universo. La solución fue la caracterización de un sistema de referencia, un campo vectorial temporal cuyas curvas integrales describen el movimiento *medio* de la materia cósmica en el Universo. Las velocidades peculiares de objetos concretos tales como galaxias pueden medirse, y despreciarse en magnitud, respecto al sistema dado. La noción de *espacio homogéneo e isotrópo* es entonces simplemente el hecho de que las hipersuperficies ortogonales al campo vectorial –cuya existencia se supone– tienen dicha propiedad. Bello y matemáticamente simple.

[11]Lanczos es muy conocido para los relativistas, colaboró con Einstein, de quien fue asistente en 1928-29 en Berlin, e hizo múltiples contribuciones a la relatividad general a lo largo de su vida.

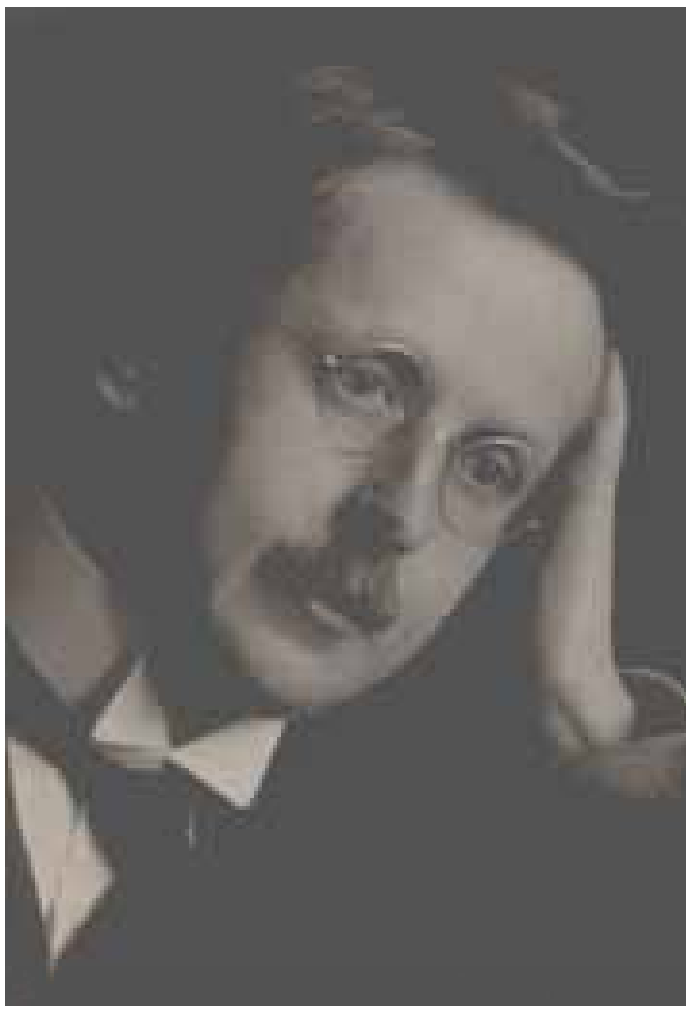
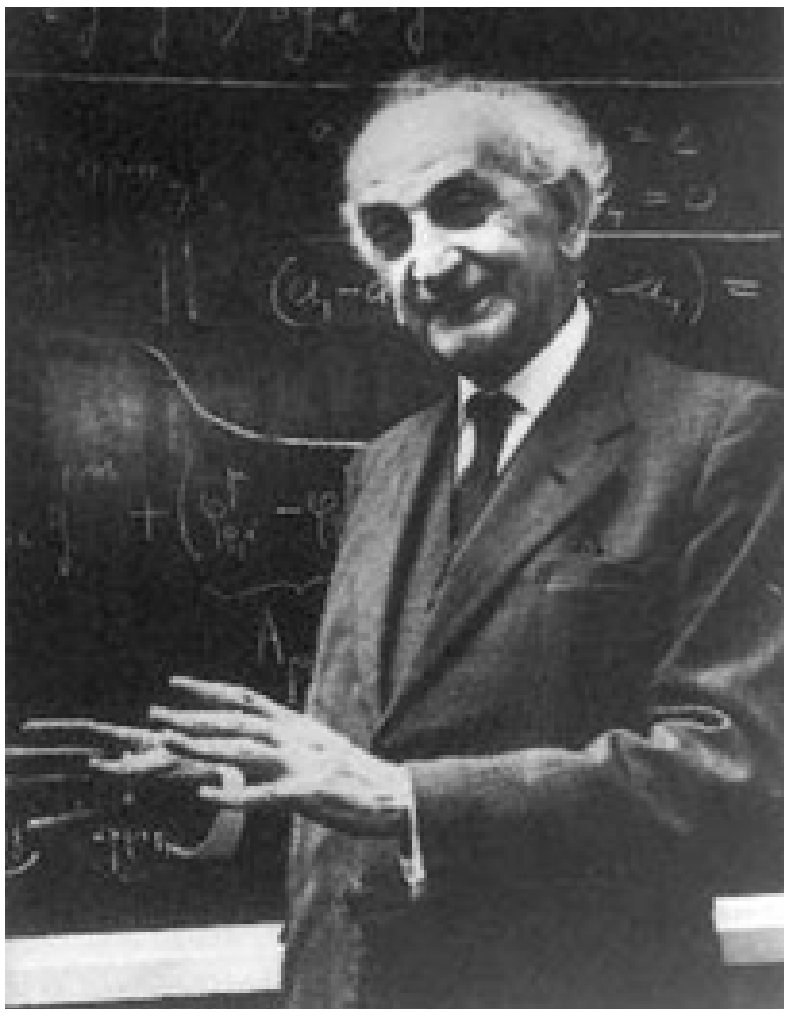

H. Weyl (izquierda) y C. Lanczos.

Con notación y visión matemáticas actuales, la cuestión es que el espacio-tiempo de de Sitter (17) es de curvatura constante positiva por lo que de hecho se puede dar en su forma canónica de Riemann, u otras. Es evidente, por otro lado, que la forma (17) no cubre toda la variedad $\mathbb{R}\times S^3$. Una forma adecuada del espacio-tiempo de curvatura constante positiva es la siguiente

$$g = -c^2 d\bar{t}\otimes d\bar{t} + a^2\cosh^2(c\bar{t}/a)\left[d\bar{\chi}\otimes d\bar{\chi} + \sin^2\bar{\chi}\left(d\theta\otimes d\theta + \sin^2\theta\, d\phi\otimes d\phi\right)\right]$$

que se puede obtener a partir de (17) mediante el cambio de coordenadas

$$t = \frac{a}{c}\log\left[\frac{\sinh(c\bar{t}/a) + \cosh(c\bar{t}/a)\cos\bar{\chi}}{\sqrt{1-\cosh^2(c\bar{t}/a)\sin^2\bar{\chi}}}\right], \qquad \sin\chi = \cosh(c\bar{t}/a)\sin\bar{\chi}$$

proporcionando de hecho una extensión máxima de (17). Como se aprecia, esta forma del espacio-tiempo de de Sitter es un caso particular también de los modelos de Friedman (19), donde $a(\bar{t}) = a\cosh(c\bar{t}/a)$.

En esta forma no estática del campo métrico la explicación dada en la sección 5.1 de la aparición del corrimiento al rojo no se aplica, ya que ahora la componente $g_{\bar{t}\bar{t}}$ es constante. Desde un punto de vista físico esto es una manifestación de cambio de sistema de observadores fundamentales, o en breve, de sistema de referencia, ya que se ha pasado de un campo vectorial temporal que respeta el Principio Cosmológico, a *otro* campo vectorial que también lo hace[12]. La explicación y fórmulas dadas por Weyl para $z$ son, esencialmente, las que se presentan en la sección 6 para los modelos corrientes en general.

### 5.5 La síntesis de Lemaître: la ciencia de la Cosmología

La síntesis definitiva que hizo encajar todas las piezas del puzzle cosmológico, combinando espléndidamente los resultados experimentales de Hubble con los modelos de Friedman, vino de la mano de George Lemaître. En un artículo esencial [15], cuyo título era meridianamente claro, certero y toda una declaración de intenciones:

*Un universe homogène de masse constante et de rayon croissant, rendant compte de la vittese radiale des nébuleuses extragalactiques*

Lemaître mejoró y generalizó los modelos de Friedman para el caso de un fluido perfecto, o sea tal que el tensor energía-momento se puede escribir como

$$\boldsymbol{T} = \varrho c^4 dt \otimes dt + p\ {}^3 g\,, \tag{21}$$

donde $p$ es la presión. El autor además hizo estimaciones del tamaño del Universo, y de su expansión, y estudió estos modelos en el contexto de las observaciones cosmológicas, dando una explicación teórica a los desplazamientos al rojo observados por Slipher y Hubble.

De nuevo los entresijos de este episodio y su posterior influencia no tienen desperdicio pedagógico. Todavía en 1927, y a pesar de los trabajos de Friedman, Lemaître, Lanczos y Weyl, además de los resultados de Slipher y Hubble, la inmensa mayoría de científicos seguía pensando que el Universo era estacionario. Por ello, los resultados de Friedman y Lemaître fueron ignorados durante años, cuando no abiertamente criticados. Por ejemplo, en el congreso Solvay de 1927 Einstein y Lemaître tuvieron oportunidad de verse y discutir, y se sabe que aquél continuaba rechazando un modelo de universo en expansión por no sustentarse en base física posible. Por comentarios de Lemaître se sabe que *grosso modo* Einstein insistía entonces en la validez matemática de los resultados, pero abominaba de la física subyacente –un eco de su frase tachada en [7, 8].

[12] De hecho, existen aún otros sistemas de observadores fundamentales en el espacio-tiempo de de Sitter. Un ejemplo de relevancia histórica y física está dado por la forma

$$g = -c^2 dT \otimes dT + e^{2cT/a}\,(dx \otimes dx + dy \otimes dy + dz \otimes dz)$$

que es la base del modelo del estado estable de Bondi, Gold y Hoyle, el mayor rival para los modelos de la gran explosión. No es éste el lugar para una discusión de estos temas. En todo caso, y naturalmente, hay que decir que no todos ellos cubren toda la variedad.

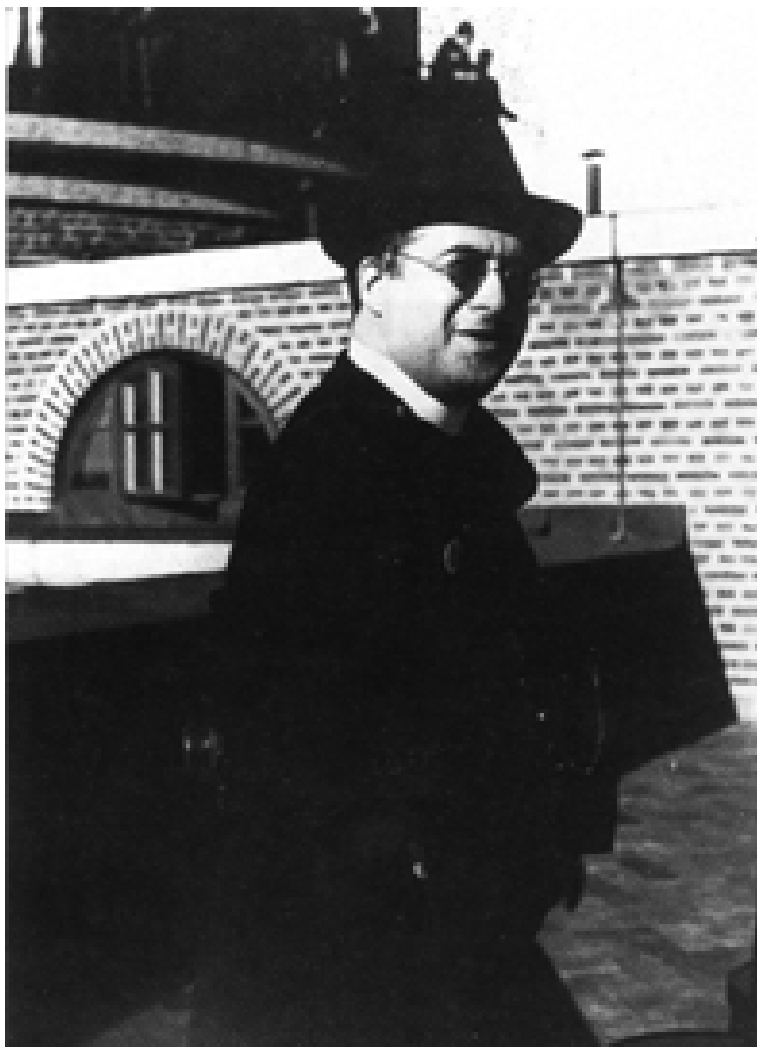

G. Lemaître

Las cosas cambiaron paulatinamente a partir de un encuentro de la Royal Astronomical Society, en Enero de 1930, cuando Eddington[13] y otros [17] empezaron a sopesar la posibilidad de soluciones no estacionarias, o sea, otras soluciones aparte de (12) y (17). Curiosamente, Lemaître había sido brevemente un *postdoc* con Eddington, y éste había visto su trabajo [15], pero *lo había olvidado por completo.* Cuando Lemaître vio las observaciones de Eddington en [17], le escribió recordando la existencia de su trabajo en [15], con el cual Eddington quedó, ahora sí, profundamente impresionado y empezó a mover los hilos para darle la publicidad merecida. De esta manera, Eddington por un lado probó la inestabilidad del universo estático de Einstein [18] y por otro arregló las cosas para que una traducción al inglés de [15] se publicase en [19].

Teniendo en cuenta que los resultados de Hubble se habían asentado, y que Lemaître pudo usar su modelo para explicarlos, por fin llegó el momento en el que la fruta estaba madura y todos finalmente aceptaron lo inevitable: *la expansión del Universo, tal y como predecía la Relatividad General.* En particular, Einstein mismo reconoció en 1931 [20] que estos modelos daban una explicación plausible de la expansión universal (o mejor, de las medidas de $z$), e incluso escribió que por ello el término cosmológico en (13) era superfluo y ya no estaba justificado. Si esto es así o no ha provocado ríos de tinta, amargas discusiones, y problemas aún hoy irresolutos.

En resumen, se había conseguido un marco teórico *falsable*, que predecía fenómenos después observados y explicaba otras observaciones previas: la Cosmologia adquiría, al fin, carácter científico.

## 6 Los modelos corrientes

Como era de esperar, hay muchas otras curiosidades, intentos fallidos, y vaivenes diversos en esta ya larga historia. Dos que merece la pena mencionar son (i) una vez que los resultados de Lemaître fueron aceptados como clarificadores y correctos, hubo una tendencia –pensemos que desinteresada y honrada– a olvidar los artículos originales de Friedman; y (ii) que ninguno de los dos, Friedman y Lemaître, encontraron el caso más sencillo de todos, aquel en que el campo métrico espacial ${}^3g$ asociado al sistema de referencia cósmico es plano.

[13]A.S. Eddington era uno de los astrónomos más influyentes y respetados de la época, aparte de un relativista de pro, siendo el director de las expediciones para observar los eclipses que confirmaron la desviación de la luz en acuerdo con la relatividad general. Además escribió uno de los primeros tratados importantes sobre relatividad general [16].

Esto viene a cuento porque la primera persona en encontrar este modelo sencillo fue el matemático H.P. Robertson [21] en 1929. El propio Robertson [22] y casi simultáneamente otro matemático, A.G. Walker [23], fueron los primeros en derivar de manera detallada y rigurosa todos los modelos posibles que satisfacen el Principio Cosmológico. Hoy en día se suele llamar a estos modelos de *Friedman-Lemaître-Robertson-Walker* (FLRW), o abreviadamente *modelos corrientes* (*standard models* en inglés). Robertson y Walker además caracterizaron las propiedades de simetría de estos modelos, así como sus posibles inmersiones en espacios de mayor dimensión.

La presentación habitual que se hace actualmente de estos modelos se basa simplemente en el Principio Cosmológico, o sea en la suposición de que existe un campo vectorial temporal $\vec{u}$ con hipersuperficies ortogonales ($g(\cdot,\vec{u}) \propto dt$) que son máximamente simétricas, o sea tal que el espacio en cada instante de tiempo propio ($t =$cte. ) del observador asociado a $\vec{u}$ es homogéneo e isótropo. Por ello, el espacio-tiempo admite un grupo de isometrías de 6 parámetros cuyas superficies de transitividad son las mencionadas hipersuperficies $t =$cte.

Bajo esta suposición, el campo tensorial métrico se puede escribir localmente en un sistema de coordenadas $\{t, \chi, \theta, \phi\}$ como

$$g = -c^2\, dt \otimes dt + a^2(t) \left[ d\chi \otimes d\chi + \Sigma^2(\chi, k) \left( d\theta \otimes d\theta + \sin^2\theta\, d\phi \otimes d\phi \right) \right] \quad (22)$$

donde $a(t)$ es una función cualquiera llamada *factor de escala*, $k = \pm 1, 0$ es un *índice de curvatura* que selecciona los tres posibles signos para la curvatura constante de las hipersuperficies $t =$cte. , y la función $\Sigma(\chi, k)$ es

$$\Sigma(\chi, k) \equiv \begin{cases} \sin\chi & \text{if} \quad k = 1 \\ \chi & \text{if} \quad k = 0 \\ \sinh\chi & \text{if} \quad k = -1. \end{cases}$$

Por razones obvias se suelen denominar modelos cerrado, abierto y espacialmente plano a los que tienen $k = 1, -1, 0$ respectivamente. Las variedades subyacentes tradicionales[14] son $\mathbb{R}^4$ para los casos $k = 0, -1$ y $\mathbb{R} \times S^3$ para $k = 1$. Es elemental comprobar que los únicos tensores energía-momento $\boldsymbol{T}$ compatibles con (22) a través de las ecuaciones de Einstein (13) son del tipo

[14]Hay que hacer constar, naturalmente, que se puede variar la topología del espacio en estos modelos sin mayores problemas, y así se puede tener un modelo 'espacialmente plano' con topología $\mathbb{R} \times S^1 \times S^1 \times S^1$, por poner sólo un ejemplo.

(21), y dichas ecuaciones conducen simplemente a

$$\frac{8\pi G}{c^2}\varrho + \Lambda = 3\frac{\dot{a}^2 + k}{a^2}, \tag{23}$$

$$\frac{8\pi G}{c^4}p - \Lambda = -2\frac{\ddot{a}}{a} - \frac{\dot{a}^2 + k}{a^2}\,. \tag{24}$$

Estas son las leyes que gobiernan la forma y evolución del Universo y su contenido material. Su extrema simplicidad permite analizar propiedades generales de los modelos corrientes con métodos elementales. Las más sobresalientes se siguen de la fórmula

$$\frac{8\pi G}{c^4}(\varrho c^2 + 3p) - 2\Lambda = -6\frac{\ddot{a}}{a}, \tag{25}$$

que es una consecuencia inmediata de (23-24). En un universo típico, formado por galaxias, sus aglomeraciones y radiación, es seguro que la cantidad $\varrho c^2 + 3p$ es positiva[15], de donde se sigue inmediatamente la siguiente conclusión sorprendente para universos en expansión ($\dot{a}_0 > 0$, donde el subíndice cero indica el valor actual): si $\Lambda \leq 0$, *existe un valor finito* $\tilde{t} < t_0$ *tal que* $a(t \to \tilde{t}) \to 0$[16]. Además, es obvio que

$$t_0 - \tilde{t} \leq \frac{1}{c}\frac{a}{\dot{a}}(t_0) \equiv \frac{1}{H_0}$$

donde la nomenclatura escogida para $H_0$, en concordancia con la notación de (20), no es causalidad como ahora veremos. Hay que notar, por lo tanto, que en los modelos corrientes la edad del Universo queda directamente relacionada con un parámetro observable y medido con gran precisión, la constante de Hubble. En todo caso, la catástrofe que aparecía en los modelos de Friedman, y discutida anteriormente, parece ser bastante genérica en los modelos corrientes con una constante cosmológica no positiva, en particular si $\Lambda = 0$.

Una pregunta que surge rápidamente es ¿cuándo ocurrirá que un modelo cambie su etapa de expansión por una de contracción? Claramente, para ello hace falta que $\dot{a} = 0$ para algún valor de $t$. Usando ahora la ley (23) y el hecho de que $\varrho > 0$ es evidente que esto puede ocurrir solamente si se da una de las dos circunstancias siguientes

- $k = 1$, o sea, el modelo es cerrado;
- $\Lambda < 0$.

[15] Esta propiedad se suele llamar *condición fuerte de la energía* y, aunque se satisface ampliamente, no está asegurada si existen campos escalares o materia más exótica.

[16] Esto es de hecho ni más ni menos que un ejemplo trivial de los famosos teoremas de singularidades.

Por ello se suele identificar la posibilidad de recolapso del Universo con que el modelo sea cerrado, lo cual claro está es incorrecto en general –se hace cuando se consideran las ecuaciones originales (11), o sea, si la constante cosmológica no se toma en consideración, lo que entonces es lícito–.

Dicho esto, es natural cuestionarse si podemos discernir cuándo un universo descrito por estos modelos será cerrado o abierto. De nuevo usando la primera ley (23) y tomándola en el tiempo actual $t_0$ se deduce que

$$\text{signo}(k) = \text{signo}\left[\varrho(t_0) - \frac{H_0^2}{8\pi G} - \frac{\Lambda c^2}{8\pi G}\right].$$

La cantidad $H_0^2/8\pi G$, que como se ve es también directamente observable al estar relacionada únicamente con el parámetro observacional $H_0$, se denomina *densidad crítica.* En el caso de no considerar la constante cosmológica, esta valor decide, por comparación con la densidad real de la materia del Universo, si éste es cerrado, plano o abierto.

Por completitud, queda finalmente por probar que la ley de Hubble (20) se deduce aproximadamente en estos modelos y que de ello se sigue la relación anunciada

$$H_0 = c\frac{\dot{a}}{a}(t_0) = \frac{1}{a(t_0)}\frac{da}{dt}(t_0)\,. \tag{26}$$

El cálculo se puede hacer escogiendo un evento cualquiera $(t_e, \chi_e, 0, 0)$ desde el que se nos envía radiación electromagnética que recibimos ahora en $(t_0, \chi_0, 0, 0)$. Teniendo en cuenta que la radiación se propaga por geodésicas luminosas, es un ejercicio sencillo probar que a lo largo de tales curvas se satisface

$$|\chi_0 - \chi_e| = c\int_{t_e}^{t_0} \frac{dt}{a(t)}\,.$$

Por consiguiente, si en un instante cercano posterior $(t_e+\Delta t_e, \chi_e, 0, 0)$ el emisor nos sigue mandando radiación, que nos llegará en $(t_0+\Delta t_0, \chi_0, 0, 0)$, la fórmula anterior para ambos casos nos lleva inmediatamente a

$$\int_{t_e}^{t_0} \frac{dt}{a(t)} = \int_{t_e+\Delta t_e}^{t_0+\Delta t_0} \frac{dt}{a(t)}\,.$$

Esta relación implica a primer orden en $\Delta t_0, \Delta t_e$

$$\frac{\Delta t_0}{\Delta t_e} = \frac{a(t_0)}{a(t_e)} = \frac{\nu_e}{\nu_0} \equiv 1 + z$$

de donde, usando la relación que se puede obtener para $t_0(t_e)$ deducida de la ecuación de las geodésicas luminosas y haciendo un desarrollo en serie se obtiene finalmente

$$z \approx \dot{a}(t_0)\,|\chi_0 - \chi_e|\,.$$

Calculando la distancia actual $D_{0e}$ entre los dos puntos con el campo métrico (22) se sigue fácilmente

$$D_{0e} = a(t_0)\,|\chi_0 - \chi_e|$$

de donde finalmente llegamos a la relación (20) con la definición (26), como se deseaba.

## 7 La irrupción de Gödel

A partir de 1933 hubo diversas generalizaciones y variados estudios de los modelos corrientes, a parte de nuevas soluciones de las ecuaciones de Einstein con simetrías esférica y cilíndrica. Siendo estos pasos de importancia, no podemos tratarlos aquí en aras de la concisión. Por ello, damos el salto a 1949 cuando apareció un número del *Reviews of Modern Physics* dedicado a Einstein por su septuagésimo cumpleaños. En este número, el mejor lógico matemático del siglo XX y amigo personal –aparte de colega en Princeton– de Einstein, Kurt Gödel (Brünn, imperio austro-húngaro, 1906 – Princeton, EE UU, 1978), publicó un artículo sobre cosmología relativista que zarandeó las creencias de los físicos relativistas, hizo zozobrar la base *causal* de la teoría, e incitó numerosos avances y nuevos conceptos fundamenales que se desarrollarían en la segunda mitad del siglo pasado.

Una vez más, vamos a ver como la falta de prejuicios físicos de un matemático, como ya ocurrió en el caso de Friedman, permite una profundización, una mejora y un salto adelante gigantesco en la física relativista. La contribución de Gödel, además de ser originalísima, tuvo una influencia enorme en la comunidad relativista y se considera el origen de los estudios modernos en teoría de causalidad, técnicas globales en variedades lorentzianas, estudios de singularidades, y muchos otros relacionados.

Gödel se planteó la posibilidad de *rotación* de la materia universal. Claro está, la rotación ha de referirse a un cierto sistema de referencia. Por ello, si se quiere definir la rotación del sistema de observadores fundamental respecto del sistema de referencia que define, hay que hablar de rotación respecto de la *brújula de la inercia*, en palabras del propio autor. Aunque ya se conocían soluciones con rotación en relatividad general (debidas a Lanczos y van Stockum), la posibilidad de materia rotante en el *contexto cosmológico* no se había considerado, especialmente debido a la manera de entender el Principio Cosmológico: dado que la homogeneidad e isotropía del espacio se refería a las hipersuperficies ortogonales al campo vectorial temporal que define el sistema de observadores fundamental, este campo es *integrable*, es decir, considerado como 1-forma es proporcional a una diferencial exacta. Esto es incompatible con la rotación intrínseca de la materia. En términos matemáticos, si $\boldsymbol{u} = g(\cdot\,,\vec{u})$ es la 1-forma asociada al campo vectorial unitario fundamental $\vec{u}$, la rotación se puede definir como

$$\boldsymbol{\omega} \equiv *(\boldsymbol{u} \wedge d\boldsymbol{u})$$

donde $*$ es el dual de Hodge. Trivialmente, $\boldsymbol{\omega}$ es espacial y ortogonal a $\vec{u}$, y como probó Gödel coincide con la velocidad angular newtoniana en coordenadas adaptadas al campo $\vec{u}$. Como es evidente, si $\boldsymbol{u}$ es integrable la rotación se anula.

El modelo de Gödel [24] es un espacio-tiempo $(\mathbb{R}^4, g)$ cuyo campo tensorial métrico se expresa en coordenadas cartesianas $\{t, x, y, z\}$ como sigue

$$g = -(cdt + e^{ax}dz) \otimes (cdt + e^{ax}dz) + dx \otimes dx + dy \otimes dy + \frac{1}{2}e^{2ax}dz \otimes dz \tag{27}$$

siendo $a$ una constante arbitraria. Obviamente, el caso $a = 0$ es el espacio-tiempo plano de Minkowski. Este campo métrico satisface las ecuaciones de Einstein (13) con el miembro derecho de un fluido sin presiones dado por

$$\boldsymbol{T} = \varrho c^2(cdt + e^{ax}dz) \otimes (cdt + e^{ax}dz)$$

donde además se satisfacen las relaciones

$$\frac{4\pi G}{c^2}\varrho = -\Lambda = g(\boldsymbol{\omega}, \boldsymbol{\omega}) > 0$$

de donde se sigue que en este caso la constante cosmológica es negativa[17].

Si uno descompone el campo tensorial $\nabla \boldsymbol{u}$ en sus partes ortogonales y tangentes a $\boldsymbol{u}$, y las primeras en su parte simétrica y antisimétrica se obtiene en general

$$\nabla \boldsymbol{u} = -\boldsymbol{u} \otimes \boldsymbol{a} + \boldsymbol{\Sigma} + *(\boldsymbol{u} \wedge \boldsymbol{\omega})$$

donde $\boldsymbol{a} = \nabla_{\vec{u}}\boldsymbol{u}$ es la aceleración (o primera curvatura) y el tensor simétrico espacial $\boldsymbol{\Sigma}$ se suele denominar tensor de deformación. En el universo (27), el campo vectorial fundamental está dado por

$$\vec{u} = \frac{1}{c}\frac{\partial}{\partial t}, \qquad \boldsymbol{u} = -(cdt + e^{ax}dz)$$

y Gödel acertadamente mencionó que tanto la aceleración como la deformación se anulan, solamente la rotación es no nula. Estas nociones, que son cantidades cinématicas, fueron posteriormente base de muchos desarrollos y ampliamanete usadas en la Cosmología teórica.

El hecho de que la deformación y la aceleración sean cero implican, en particular, que $\vec{u}$ es un campo de Killing[18]. Por ello, el modelo de Gödel es

[17] El universo de Gödel también se puede interpertar como solución para un fluido perfecto.

[18] Un campo de Killing es el generador local de un grupo uniparamétrico de *isometrías* locales.

estacionario. Se puede probar además que los campos vectoriales siguientes

$$\frac{\partial}{\partial x} - az\frac{\partial}{\partial z}, \qquad -2e^{-ax}\frac{1}{c}\frac{\partial}{\partial t} + az\frac{\partial}{\partial x} + \left(e^{-2ax} - \frac{a^2z^2}{2}\right)\frac{\partial}{\partial z}, \qquad \frac{\partial}{\partial y}, \qquad \frac{\partial}{\partial z}$$

son también campos de Killing, todos ellos linealmente independientes. Esto implica, como bien señaló Gödel, que el espacio-tiempo es homogéneo, con un grupo de isometrías de 5 parámetros que actúa de manera multiplemente transitiva.

Todas las propiedades anteriores eran relativamente novedosas, especialmente por su presentación, pero la propiedad auténticamente sorprendente, demoledoramente impactante, que probó Gödel es que el espacio-tiempo (27) es causalmente orientable –o sea, se puede asociar consistentemente una dirección al futuro, de manera continua, en toda la variedad; todos los vectores tangentes temporales o luminosos se pueden dividir en futuros y pasados, consistentemente de manera que las curvas temporales quedan también así orientadas, y de forma continua–, y a pesar de eso existen *curvas temporales cerradas dirigidas al futuro*. Hablando llanamente, dejando pasar el tiempo como es natural, ¡uno puede llegar a su propio pasado! Más aún, diseñando adecuadamente el camino del viaje, uno puede llegar tan atrás al pasado de sí mismo –yendo al futuro continuamente– como se quiera.

Varias paradojas surgen inmediatamente, como la tradicional de "persona que va al pasado, mata a su madre antes de que mantenga la relación sexual correspondiente a su engendramiento, y... ¿?". En fin, esta posibilidad, en una solución de las ecuaciones que además está libre de singularidades, es estacionaria y con un contenido material razonable, es muy desagradable, molesta y perturbadora. Supongo que esto es así incluso para muchos matemáticos, por lo que me permito mostrar explícitamente la siguiente familia de curvas

$$ct = A\left(2\sin\tau - \sin\tau\cos\tau\right), \qquad x = -B\cos\tau, \qquad y = 0, \qquad z = -2A\sin\tau$$

donde $\tau$ es el parámetro de la curva, y dejo como ejercicio al lector que compruebe para qué valores de las constantes $A$ y $B$ son temporales, futuras y cerradas –topológicamente $S^1$– en el universo de Gödel dado por (27).

Hay que resaltar, para mayor desconsuelo, que dado que (27) es homogéneo, existen curvas temporales futuras cerradas pasando por cualquier punto de la variedad. Y también que no aparecen debido a raras propiedades topológicas, ni se pueden evitar pasando a supuestos "recubridores universales", porque la variedad es en este caso $\mathbb{R}^4$, que es simplemente conexa.

Gödel relacionó estas circunferencias temporales con la ausencia de hipersuperficies ortogonales a $\boldsymbol{u}$, o sea, con la rotación, y claramente demostró que no podía haber una función temporal, o sea, una función continua sobre $\mathbb{R}^4$ que crezca en la dirección futura a lo largo de *cualquier* curva temporal futura[19]. Es evidente que esta ausencia es una condición necesaria para la existencia de curvas temporales cerradas, y por ello esta discusión fue la base de gran parte de las posteriores condiciones de causalidad que los físicos trataron de definir para evitar las paradojas mencionadas. Esta teoría de la causalidad se reveló como de fundamental importancia para los desarrollos relacionados con los teoremas de singularidades, o la (in)completitud geodésica de los espacio-tiempos. Notemos, de paso, que la solución (27) es geodésicamente completa.

Por si todo lo anterior no fuera suficiente, Gödel ya menciona en su artículo los resultados de Bianchi sobre espacios tridimensionales que admiten grupos continuos de transformaciones. Dice por ejemplo que la parte de (27) en la que se prescinde del término $dy \otimes dy$ es esencialmente uno de los espacios en forma canónica encontrada por Bianchi en [26]. El "esencialmente", claro está, se refiere a la diferente signatura.

Esta línea de actuación la generalizó el propio Gödel en su segundo artículo sobre relatividad [25], otro auténtico punto de inflexión, conteniendo un marasmo de ideas y de nociones que fueron de utilidad e inspiración para los físicos durante décadas posteriores. Habida cuenta de que su modelo (27) no era realista al no expandirse (es estacionario), trató de encontrar nuevas soluciones que se expandieran (lo que necesita que la traza de la deformación $\mathbf{\Sigma}$ sea no nula) y mantuvieran la rotación. Para ello consideró modelos que fueran invariantes por el grupo $SO(3)$, actuando de manera transitiva sobre hipersuperficies espaciales.

Esta es la base, cambiando simplemente el grupo, de los modelos espacialmente *homogéneos* (no necesariamente isótropos) que se han venido en llamar *modelos de Bianchi* desde la década de los 1950. Luigi Bianchi [26] clasificó, a finales del siglo XIX, todos los espacios riemannianos tridimensionales que admitían un grupo de isometrías continuo. En particular, si la acción es transitiva, por lo que las órbitas son tridimensionales, obtuvo la forma canónica en coordenadas adaptadas a los campos de Killing de estos espacios, que resulta deducirse de la clasificación de los propios grupos de 3 parámetros, de acuerdo con sus constantes de estructura. Estos grupos se denominan tipos de Bianchi y se clasifican desde el I hasta el IX, aunque alguno de los casos son en realidad familias uniparamétricas de constantes de estructura. Esta misma nomenclatura se heredó en relatividad y los modelos cosmológicos de Bianchi se dividen en grupos del tipo I al IX.

La idea es bastante simple. Dado el grupo triparamétrico actuando en el espacio-tiempo con órbitas tridimensionales espaciales, y sus constantes de estructura, se puede probar que las 1-formas invariantes por la izquierda (o la derecha) del grupo se trasladan a la variedad, mediante la acción del grupo, a 1-formas $\{\boldsymbol{\theta}^\alpha\}$ $(\alpha, \beta = 1, 2, 3)$ que son invariantes por los correspondientes campos de Killing que generan la acción del grupo. Es fácil demostrar entonces que el campo métrico adopta necesariamente la forma

$$g = -c^2 dt \otimes dt + g_{\alpha\beta}(t)\boldsymbol{\theta}^\alpha \boldsymbol{\theta}^\beta \tag{28}$$

donde las funciones $g_{\alpha\beta}$ dependen ahora de $t$. La forma concreta de las 1-formas $\{\boldsymbol{\theta}^\alpha\}$ depende del grupo escogido y su expresión explícita es sencilla [26], véase, por ejemplo, [27]. Es obvio que $\boldsymbol{u} = -cdt$ define un sistema de referencia cuyos espacios ortogonales son homogéneos. Pero eso, como se ha explicado antes, impide que una hipotética materia moviéndose a lo largo de las curvas integrales de $\vec{u}$ tenga rotación. Por ello, aunque lo lógico parece ser escoger en estos modelos $\boldsymbol{u} = -cdt$ como la 1-forma media del movimiento de la materia, Gödel demostró la existencia de soluciones para un fluido perfecto donde el campo de velocidades es distinto de $\vec{u}$. La distribución material tiene rotación, y este tipo de modelos, los de mayor complejidad dentro de la familia de Bianchi, se denominan ahora *modelos ladeados* ("tilted models" en inglés).

[19]De paso, Gödel también probó que no existe ninguna hipersuperficie sin borde espacial en (27), otra propiedad inesperada.

Los modelos de Gödel en [25] eran del tipo Bianchi IX, y escogió ese grupo porque deseaba que las órbitas, o sea, el espacio, fuera compacto. Hizo notar, cosa que es válida para todos los modelos de Bianchi, que las ecuaciones de Einstein se reducen a un sistema diferencial *ordinario* para las variables $g_{\alpha\beta}(t)$, e incluso dio un lagrangiano para derivar estas ecuaciones en su caso. Obviamente, esta propiedad de los modelos de Bianchi tiene una importancia enorme, y permite utilizar todas las técnicas modernas de sistemas dinámicos, e incluso las más modernas teorías de caos, aplicadas a modelos de Universo. Se puede ver un resumen relativamente reciente en [28].

El segundo artículo de Gödel está escrito con un estilo un tanto misterioso, enumerando gran cantidad de resultados y propiedades pero sin pruebas o con levísimas indicaciones de las mismas. Esto conformó un rompecabezas que puso manos a la obra a grandes físicos y algunos matemáticos, que permitieron probar los antedichos resultados y desarrollar las brillantes ideas godelianas. Entre ellos hay que citar a A. Taub, quien casi simultáneamente a [25] encontró todas las formas canónicas de los modelos de Bianchi [29], estudió cuándo eran soluciones sin materia y dio una detallada e iluminadora discusión de sus propiedades físicas.

Un resumen de las ideas germinales de Gödel, que se convertirían en aspectos esenciales de la relatividad moderna y de los desarrollos actuales, es (véase [30] para más detalles):

- La rotación de la materia que define el movimiento medio universal respecto de la brújula de la inercia. Esto conduce al abandono definitivo de las ideas machianas.

- Posibilidad de curvas temporales cerradas. Relación de esto con la ausencia de funciones temporales globales. Esto reveló la necesidad de una teoría de la causalidad, que fue posteriormente desarrollada y es base de la geometría lorentziana global.

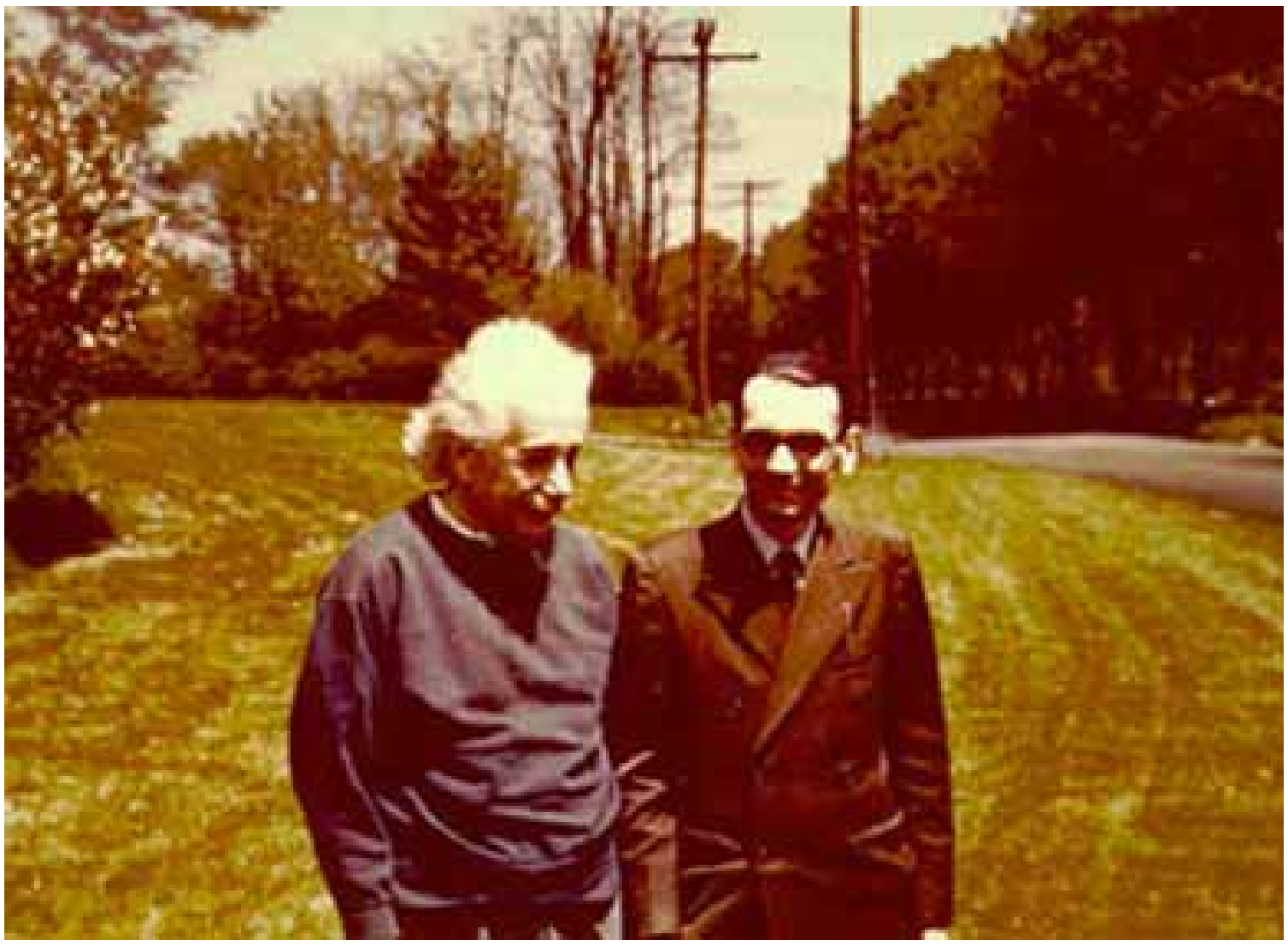

A. Einstein y K. Gödel en Princeton.

- Uso de la clasificación de Bianchi, y de las acciones de grupos sobre el espacio-tiempo. Esto, junto con la clasificación algebraica del tensor de curvatura, despejó el camino para la clasificación y proliferación de soluciones de las ecuaciones de Einstein.
- Modelos ladeados, y uso por tanto de coordenadas no adaptadas al fluido.
- Definición de las cantidades cinemáticas, lo que parece ser tuvo una influencia determinante en la aparición de la ecuación de Raychaudhuri y la correspondiente focalización de geodésicas.
- El punto anterior es la base, junto con los dominios de dependencia, también subyacentes en el trabajo de Gödel, de los famosos teoremas de singularidades.

Todos estas líneas fueron la base del gran desarrollo de la relatividad moderna –debida entre otros a Penrose, Hawking, Geroch, Carter–, en la segunda mitad del siglo XX.

Pero esa es otra historia...

## 8 Epílogo

La principal conclusión de esta larga narración es, en mi modesta opinión, que no tendríamos una ciencia de la Cosmología si Einstein no hubiera desarrollado su teoría de la relatividad general, aplicándola a todo el Universo en 1917 [1]. Siendo las observaciones de Slipher y Hubble absolutamente esenciales, era por otro lado imprescindible tener una teoría que ligase la forma geométrica del Universo con su contenido, distribución material, y evolución temporal. En consecuencia, la cosmología se puede apuntar al ya largo haber de las revoluciones científicas de su responsabilidad.

Otra conclusión de interés es que la aportación de los matemáticos fue determinante, aparte de muy original, para llevar a buen puerto los incipientes primeros pasos teóricos y para una mejor comprensión de los resultados obtenidos. La falta de prejuicios físicos, como ya he mencionado varias veces, permitió a grandes matemáticos realizar los avances que se les escapaban a los desconcertados físicos teóricos.

Lamentablemente, la aportación de los matemáticos a la relatividad en la segunda mitad del siglo XX, a pesar de que gran parte de los desarrollos teóricos son fundamentales y básicamente matemáticos, no estuvo a la misma altura. Estoy seguro de que la visión original de los matemáticos, con puntos de vista diferentes, sus ángulos y perspectivas de análisis globales, y su falta de prejuicios, podría llevar a un desarrollo sin igual, a nuevos puntos de inflexión de la teoría y su sustrato geométrico, y a aportar novedosas y fructíferas ideas.

¡Ánimo!

## Referencias

José M M Senovilla
Física Teórica
Universidad del País Vasco
Apartado 644
48080 Bilbao
Correo electrónico: `josemm.senovilla@ehu.es`